\documentclass[12pt]{article}
\usepackage[utf8]{inputenc}

\usepackage{pgf,tikz}
\usetikzlibrary{decorations.pathreplacing,shapes.misc,calc}
\usetikzlibrary{arrows.meta,decorations.markings}

\usepackage{nicefrac}
\usepackage{graphicx}
\usepackage{subfigure}
\usepackage{hyperref} 
\usepackage{amsmath}
\usepackage{amssymb}
\usepackage{amsthm}
\usepackage{amsfonts} 
\usepackage{bm}
\usepackage{cases}
\usepackage{xcolor}
\usepackage{color,soul}
\usepackage{caption} 
\usepackage{graphicx}
\usepackage{epstopdf}
\usepackage{enumerate}
\usepackage{ifthen}
\usepackage{mathtools}  
\usepackage{notes2bib}
\usepackage{placeins} 
\usepackage{setspace}
\usepackage{sidecap} 
\usepackage{url}
\usepackage{ucs} 
 
\usepackage{float}

\usepackage{lineno}

\def\epsilon{\varepsilon}
\newcommand{\BF}{\bf\boldmath }

\newcommand{\ie}{{\it{i.e.}}}

\newcommand{\eps}{\varepsilon}
\newcommand{\Nbhd}{\text{Nbhd}}

\newcommand{\R}{{\mathbb{R}}}

\def\cM{{\cal M}}
\def\cN{{\cal N}}

\def\bR{{\mathbb R}}

\newcommand{\beq}{\begin{linenomath}\begin{equation*}} 
\newcommand{\eeq}{\end{equation*}\end{linenomath}} 
\newcommand{\beqn}{\begin{linenomath}\begin{equation}} 
\newcommand{\eeqn}{\end{equation}\end{linenomath}} 
\usepackage[normalem]{ulem}

\newtheorem{theorem}{Theorem}
\newtheorem{proposition}{Proposition}

\newtheorem{lemma}[proposition]{Lemma}
\newtheorem{conjecture}[proposition]{Conjecture}
\newtheorem{corollary}[proposition]{Corollary}

\def\Fto{{\mapsto}}

\definecolor{dgn}{rgb}{0.2, .6, 0.0}

\newcommand{\black}{\color{black}{}}

\newcommand{\blue}{\color{black}{}}

\definecolor{orange}{rgb}{0.8, .2, 0.0}

\definecolor{darkgreen}{rgb}{0, .6, 0.0}

\definecolor{jim}{rgb}{0,.5,.5}
\definecolor{rob}{rgb}{0,0.4, 0.6}

\definecolor{purple}{rgb}{0.1, .4, 0.0}

\newcommand{\bS}{{\mathbb{S}}}

\newcommand{\bydef}{\overset{\mathrm{def}}{=\joinrel=}}

\usepackage{enumitem}

\usepackage{authblk}

\begin{document}

\title{What is the graph of a dynamical system?}

\author[$1$]{Chirag Adwani}
\author[$2$]{Roberto De Leo}
\author[$3$]{James A. Yorke}

\affil[$^1$]{Department of Mathematics, University of Maryland College Park, MD 20742, USA, \texttt{cadwani@umd.edu}}
\affil[$^2$]{Department of Mathematics, Howard University, Washington DC 20059, USA, \texttt{roberto.deleo@howard.edu}}
\affil[$^3$]{Institute for Physical Science and Technology and the Departments of Mathematics and Physics, University of Maryland College Park, MD 20742, USA, \texttt{yorke@umd.edu}}


\maketitle

\centerline{\em In memory of Joseph Auslander (1930-2025).}

\abstract{Some of the basic properties of any dynamical system can be summarized by a graph. The dynamical systems in our theory run from maps like the logistic map to ordinary differential equations to dissipative partial differential equations.
Our goal has been to define a meaningful concept of graph of any dynamical system. 
As a result, we base our definition of ``chain graph'' on ``$\eps$-chains'', defining both nodes and edges of the graph in terms of chains.
In particular, nodes are often maximal limit sets and there is an edge between two nodes if there is a trajectory whose forward limit set is in one node and its backward limit set is in the other. 
Our initial goal was to prove that every ``chain graph'' of a dynamical system is, in some sense, connected, and we prove connectedness under mild hypotheses.}
  

\maketitle

\section{Introduction}

Graphs are often used to give an overview of a dynamical system or a semi-flow.
When we draw pictures showing fixed point repellors, saddles, and attractors, showing how their stable and unstable manifolds are related, we are showing a graph of the dynamical system. 

Our main results are the following two results, which we state so that the reader knows the direction we are headed.
The terms will be defined in Section~\ref{sec: chain graph}.
Proposition 1 is quite similar to results in literature such as Theorem~10.5 in~\cite{Rob01}.

\medskip\noindent
\begin{proposition}[\bf Existence and connectedness of the global attractor] 
\label{p1}
Assume a dynamical system has a connected trapping region $Q$.
Then the global attractor exists, is connected, and is contained in $Q$.
\end{proposition}
\medskip

\begin{theorem}[\bf Connectedness of chain graphs] 
    \label{Connectedness of graph}
    Assume a dynamical system has a global attractor that is connected. 
    Then its chain graph is connected.
\end{theorem}
\vskip 0.3cm
\noindent
{\bf Finding graphs of dynamical systems.} In~\cite{DY25}, we surveyed different ways of defining graphs of dynamical systems. 
The aim of the present article is to provide a simple intuitive type of graph (we call it ``chain graph'') that can be used effectively in applications.
We have chosen the definition of chain graphs for its simplicity and for its applicability to all dynamical systems.
We encourage researchers to publish graphs of their favorite dynamical systems, based on either numerical evidence or rigorous derivations.
The graph would reveal what the investigators have found, and might admittedly be incomplete. 
It is often impossible to be sure that one has a complete list of attractors, let alone of saddles and edges between nodes.
We have no systematic methods for identifying the graph of all dynamical systems.
Hopefully, our examples will provide investigators with some methods useful for 
identifying graphs. 
Investigators in partial differential equations have been most successful in describing complete graphs.
We hope that examples in finite dimensions will become more common.
We, in fact, were introduced to this field in our study of the logistic map. 

\medskip\noindent
{\bf Examples.} Figure~\ref{fig:cr} shows the graph of a flow on a disc.
This flow has two periodic orbits (labeled by $N_1$ and $N_2$ in figure) and a fixed point (labeled by $N_3$).
Every other orbit spirals either from $N_1$ to $N_2$ or from $N_3$ to $N_2$. 
The graph to the right of the flow portrait summarizes its dynamics in the following sense: the edges from $N_1$ and $N_3$ to $N_2$ represent the fact that, arbitrarily close to $N_1$ and $N_3$, there are orbits that asymptote to $N_2$.
Since there are no other nodes, this fully describes the qualitative dynamics of the system: each point of the system either is in a node, and so has some kind of recurrent dynamics, or goes from a node to another node, so it has some kind of gradient dynamics.


Figure~\ref{fig:cr}(right) shows the graph of the Lorenz system with standard parameter values.

Figure~\ref{fig:fiedler graphs} shows several other graphs of semi-flows. 
The rightmost one is the graph of the logistic map, which we discuss in some detail in Section~\ref{sec: logistic map}.
All other graphs, adapted from figures by Fiedler and Rocha~\cite{FR00}, encode the qualitative dynamics of some reaction-diffusion Partial Differential Equation (PDE), including the Chafee-Infante PDE we discuss in some detail in Section~\ref{sec: PDE}.
People unfamiliar with PDEs can skip those examples. 

\medskip\noindent
{\bf Overview.} The structure of the article is the following. In Section~\ref{sec: chain graph}, we introduce the chain graph, and we go over, without proofs, some of its most general properties. 
In Sections~\ref{sec: logistic map} through~\ref{sec: PDE}, we illustrate the chain graph of several important examples of dynamical systems. 
Finally, in Section~\ref{sec: proof} we present the proofs of the results presented in Section~\ref{sec: chain graph} \blue and in Section~\ref{sec: conj} we briefly present two conjectures. 

\black
\section{The chain graph of a Dynamical System}
\label{sec: chain graph}

The main object of the study of this article is 
a {\BF semi-flow $F$ on a metric space $(X,d)$}, namely a 1-parameter family of continuous maps $F^t:X\to X$ such that 
\beqn
\label{semi-flow}
F^0(x)=x\text{ and }F^t(F^s(x))=F^{t+s}(x)\text{ for all } t,s\ge0,x\in X.
\eeqn
A trajectory starting at a point $x\in X$ is at $F^t(x)$ after some time $t\ge 0$. 
We include both discrete-time maps and continuous-time processes, so $t\ge0$ means either integers or real numbers, depending on the choice of semi-flow. 

 Throughout this article, 
 when we write ``dynamical system'', the reader should assume that we are referring to a semi-flow $F$ on $X$ unless otherwise specified.
 Discussions here will usually focus on discrete-time semi-flows, and we leave it to the reader to adapt the notation to the continuous-time case.
 We use $x$, $y$ and $z$ to denote points in $X$. 
 We denote the distance between $x$ and $y$ in $X$ by $d(x,y)$ and
 we set $d(x,A)=\inf_{a\in A}d(x,a)$, $d(K,A)=\sup_{x\in K}d(x,A)$
 and $\Nbhd(Q,\delta)=\{x: d(x,Q)<\delta\}$.

\begin{figure}
 \centering
 \includegraphics[width=6.5cm]{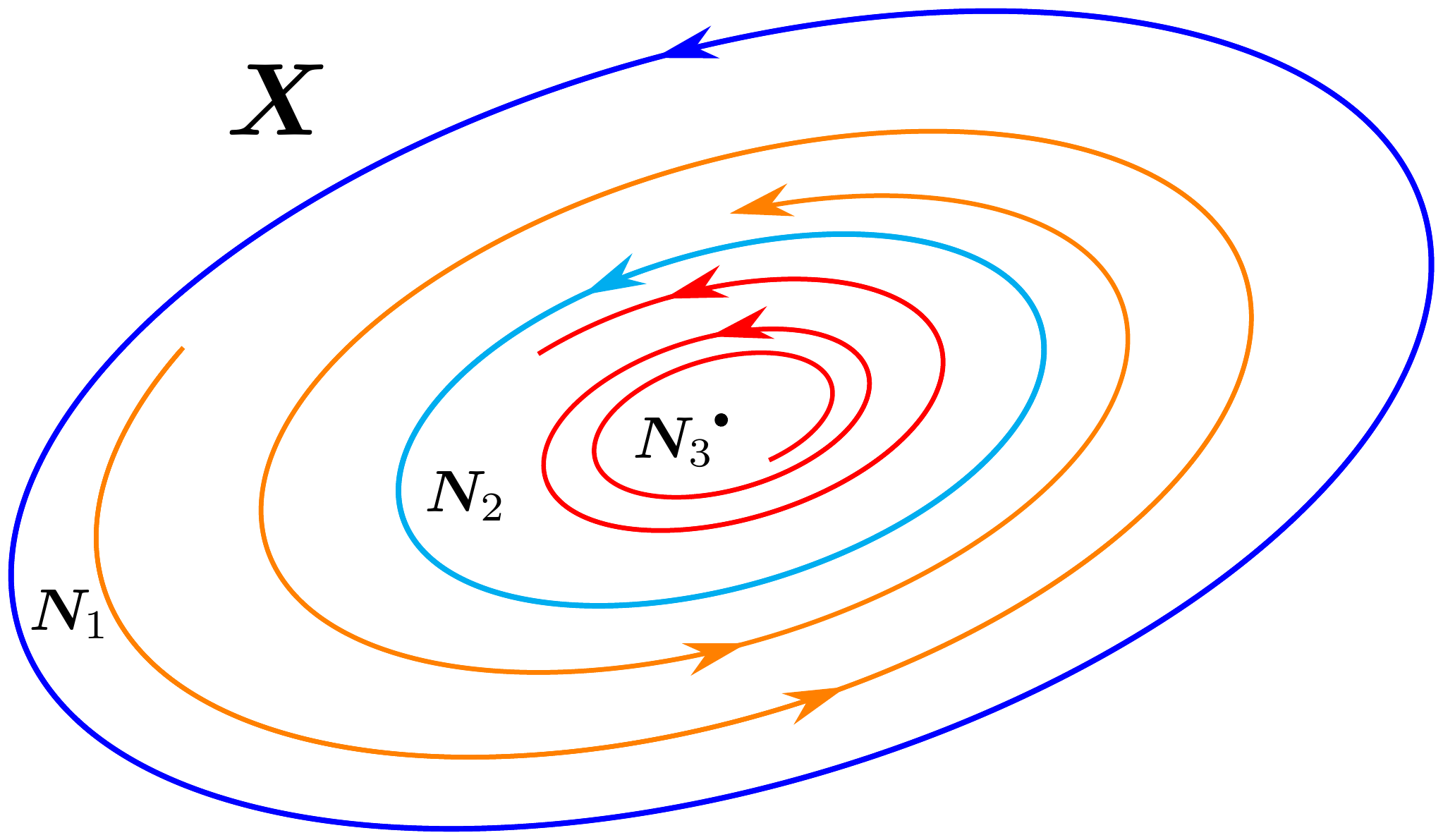}\hspace{0.3cm}
 \includegraphics[width=5.cm]{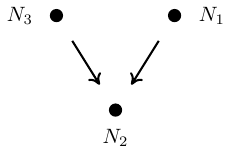}
 \caption{{\bf Representing a dynamical system by a graph. 
  (LEFT: A flow.)} Here is a dynamical system on a disk $X$. It has a fixed point $N_3$ and two periodic orbits, one of which, $N_2$, is attracting. The set $X$ is the disc bounded by the outer periodic orbit $N_1$. 
  Three nodes are visible in the picture: the outer periodic orbit $N_1$ (in blue), the inner periodic orbit $N_2$ (in cyan) and a fixed point $N_3$ (in black). 
  {\bf (RIGHT): The graph representing the flow.
 } The edges of this graph go from the repellors $N_3$ and $N_1$ to the attractor $N_2$. 
 This is also the graph of the Lorenz system with standard parameter values.
 In that case, $N_1$ and $N_3$ are two repelling fixed points and $N_2$ is the butterfly attractor (see Section~\ref{sec:Lorenz} for other Lorenz parameters).
 }
 \label{fig:cr}
\end{figure} 

 Our results hold for a wide range of dynamical systems including \black maps (e.g. see Section~\ref{sec: logistic map}), ordinary differential equations (e.g. see Sections~\ref{sec:Lorenz} and~\ref{sec: pendulum}) and reaction-diffusion partial differential equations (PDEs) such as the quasilinear equation $u_t=a u_{xx}+ f(x,u,u_x)$.
 For sake of simplicity, in Section~\ref{sec: PDE} we discuss in some detail the Chafee-Infante PDE~\cite{CI74}, which does not contain a $u_x$ term. 
 We refer the reader to Lappicy's article~\cite{Lap18} for a thorough discussion of the fully quasilinear case, including a quasilinear version of the Chafee-Infante, showing that our assumptions also hold for a wide class of quasilinear PDEs.

\medskip\noindent
{\bf Limit sets.}
The {\bf forward limit set} of $x$ under $F$ is the set
\beq 
\omega(x)\bydef\{y:\text{ there is a sequence }t_n\to\infty\text{ such that } F^{t_n}x\to y\}.
\eeq 
We also define the limit set of a set $B$ under $F$ as the set 
\beq 
\omega(B)\bydef\{y:\text{ there are sequences }t_n\to\infty, b_n\in B,\text{ such that } F^{t_n}(b_n)\to y\}.
\eeq 
When $B$ is closed and forward invariant, $\omega(B)=\displaystyle\bigcap_{t\geq0}F^t(B)$.

Notice that $\omega(x)=\omega(\{x\})$ and that $\omega(B)$ can be strictly larger than $\cup_{b\in B}\omega(b)$.
This is the case, for instance, of the ODE $x'=x(1-x)$ on $X=[0,1]$, since $\omega(X)=X$ but each point of $X$ has its forward limit set either equal to 0 or to 1.

For semi-flows, backward trajectories might not exist or they might not be unique.
We say that $\tau=(x_n)$ is a {\BF two-sided trajectory} through $x$ if the $x_n$ are defined for every $n=0,\pm1,\pm2,\dots$, $x_0=x$ and $x_{n+1}=F^1(x_n)$ for all $n$. 
Given a two-sided trajectory $\tau$ through $x$, we define the forward limit set of $\tau$ to be the forward limit set of $x$ and the {\bf backward limit set} of $\tau$ as the set 
\beq 
\alpha(\tau)\bydef\{y:\text{ there is a sequence }n_k\to-\infty\text{ such that } x_{n_k}\to y\}.
\eeq 
In case $F$ is a flow, so that through every $x$ there is a unique backward trajectory, we sometimes write $\alpha(x)$ to denote the backward limit set of the unique two-sided trajectory passing through $x$.

\medskip\noindent
{\bf Invariant sets.} 
We say that a closed set $A\subset X$ is {\bf forward invariant} (respectively, {\bf backward invariant}) under $F$ if $F^t(A)\subset A$ (respectively, $F^t(A)\supset A$) for all $t\geq0$.
We say that $A$ is {\bf invariant} if $A$ is both forward and backward invariant, namely
if $F^t(A)=A$ for all $t\geq0$.
Closed invariant sets have two properties that are quite important to us.
One is that they contain the forward limit sets of all of their points, since no point can get out of them under $F$.
The other is that each point in an invariant set $A$ has a two-sided trajectory that passes through it, since, by the definition of an invariant set, for every $x\in A$ there is a $y\in A$ such that $y=F(x)$.
Moreover, the backward limit set of each two-sided trajectory fully contained in $A$ is also fully contained in $A$ as well.

\medskip\noindent
{\bf Attracting.} 
We say that the set $A$ {\bf attracts} the set $K$ if $d(F^t(K),A)\to0$ as $t\to\infty$.
This is equivalent to the fact that, for every $\eps>0$, there is a $\tau>0$ such that $F^t(K)\subset\Nbhd(A,\eps)$ for all $t\geq\tau$.

\medskip\noindent
{\bf Global attractor.}
We say $G\subset X$ is a {\bf global attractor} if it is a maximal compact invariant set that attracts every compact set in $X$.

\smallskip\noindent
In particular, therefore, $\omega(K)\subset G$ for every compact set $K\subset X$.
Notice that such $G$, when it exists,  is unique.
Indeed, if $G'$ is another global attractor, then $G=\omega(G)\subset G'$ and $G'=\omega(G')\subset G$, so that $G=G'$. 
Similar arguments show that the global attractor contains every compact invariant set of $F$ and is contained in every set that attracts all compact subsets of $X$; see Proposition~\ref{prop: glob} in Section~\ref{sec: proof}.

\medskip\noindent
{\bf Absorbing.} 
We say that the set $A$ {\bf absorbs} the set $K$ if there is a $\tau\geq0$ such that $F^t(K)\subset A$ for every $t\geq\tau$.
\black

\medskip\noindent
{\bf Trapping regions.}
The global attractor is unique but often very hard to describe. 
We introduce here the concept of a ``trapping region'' as a tool to detect the existence of a global attractor.

\medskip
We say that a compact forward-invariant set $Q\subset X$ is a {\bf trapping region} if it absorbs each compact set in $X$.

\medskip
{\BF Throughout the article we make no assumptions on $X$
and we rather assume that $F$ has a trapping region} (which, by definition, is compact).
The assumption of having a trapping region may seem rather strong to the reader, but it is justified by the existence of many important examples, that include dissipative partial differential equations.
Hence, the topological properties of $X$ will play no role in our proofs and results.
\medskip

Proposition 1, that we already stated in the introduction and we repeat below rephrasing it slightly, shows that the existence of a trapping region implies the existence of a global attractor.

\medskip\noindent
{\bf Proposition 1 (Existence and connectedness of the global attractor).}{\em\  
Assume that $F$ has a connected trapping region $Q$.
Then the global attractor exists, is connected and equals $\omega(Q)$.
}

\medskip\black
While global attractors are unique and hard to describe, trapping regions, when they exist, abound and often it is easy to prove their existence.
If $X$ is compact, then $X$ itself is trivially a trapping region.

\smallskip\noindent
{\bf Example: the Lorenz ``butterfly'' system.}
For standard parameter values, the Lorenz ``butterfly'' system has a global attractor that includes the standard Lorenz attractor, two fixed points and their unstable manifolds. 
It is certainly complicated. 
It also has measure zero.
As Lorenz shows, though, this system has a large ball that is a trapping region~\cite{Lor63}.

\smallskip\noindent
{\bf Example: the logistic map.} 
The logistic map $\ell_a(x)=ax(1-x)$ on the non-compact set $X=(0,1)$ has a global attractor that is an interval for $1< a\leq4$.
Every trajectory of $\ell_a$ through any point in $X$ enters, in finite time, the closed invariant interval $Q_a\bydef[\ell_a^2(1/2),\ell_a(1/2)]=[\frac{a^2}{4}(1-\frac{a}{4}),\frac{a}{4}]$.
Hence, $Q_a$ is a trapping region for $\ell_a$.
Notice that here $Q_a=\omega(Q_a)$ and so $Q_a$ itself is the global attractor, showing that in some cases, global attractors can absorb (rather than just attract) each compact set of the system.

\smallskip\noindent
{\BF Example: the ODE $x'=x(1-x^2)$ on $X=\bR$.} 
Outside of the interval $[-1,1]$, each solution $x(t)$ of this ODE moves towards 0 as $t$ increases.
Hence, any interval $[a,b]$ containing $[-1,1]$ in its interior is a trapping region. 
Moreover, the reader can verify that $\omega([a,b])=[-1,1]$ if $a<-1$ and $b>1$.
Hence, $G=[-1,1]$.

\smallskip\noindent
{\BF Example: the Chafee-Infante PDE.} This PDE, investigated by Chafee and Infante in~\cite{CI74}, $u_t=u_{xx}+\lambda u(1-u^2)$, with $X=L^2(\bS^1)$, the space of integrable functions on the circle.
As discussed in Section~\ref{sec: PDE}, if $0<\lambda<1$, the global attractor of this PDE is the subset of $X$ of all constant functions $u(t,x)=c$, with $c\in[-1,1]$.

In infinite dimensions, a vast literature (e.g. see~\cite{Lad22,Rob01,Hal10,CLR12,Lap18,Lap23}) is dedicated to bounded dissipative dynamical systems, namely dynamical systems $F$ for which there exists a compact set $Q$ that absorbs every bounded set.
The set $Q$ can be assumed, without loss of generality, to be forward invariant (see~\cite{Rob01}, Exercise 10.2).
Hence, for each bounded dissipative dynamical system, if the absorbing set $Q$ is forward-invariant, then it is a global trapping region (see~\cite{Rob01}, Theorem 10.5),
and the global attractor defined in the literature above is $\omega(Q)$.

\begin{figure}
 \centering
 \includegraphics[width=13cm]{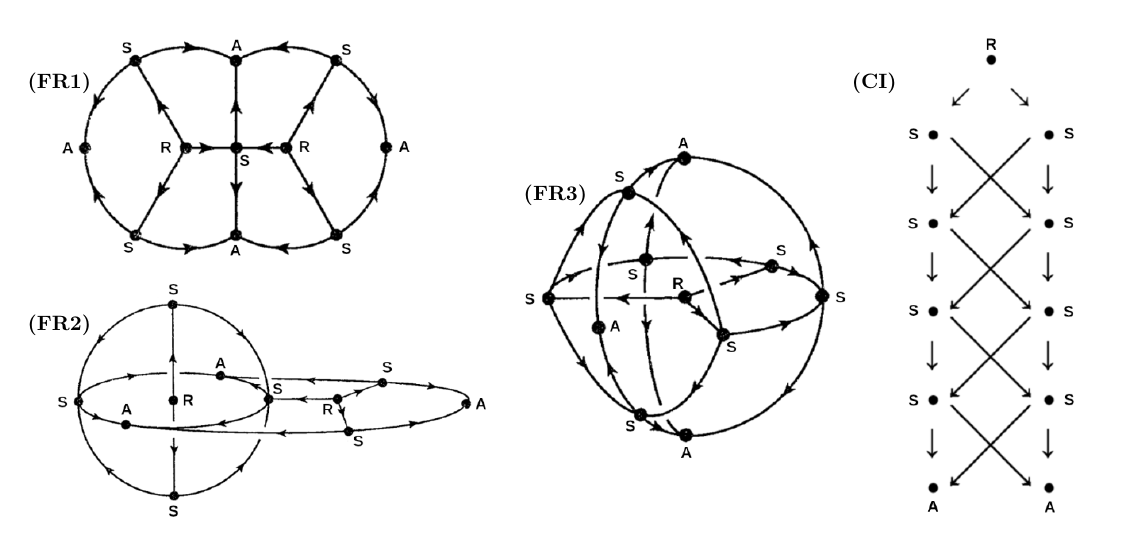}
 \caption{
 {\bf Typical figures of graphs of dynamical systems.} 
 In the graphs (FR1-FR3) and (CI), adapted from figures in Fiedler and Rocha~\cite{FR96}, each large dot represents a node, and each node is a fixed point (a steady state). 
 There is an edge from a node $N_1$ to node $N_2$ if there is a trajectory whose forward limit set is in $N_2$ and whose backward limit set is in $N_1$.
 Each of their dynamical systems is a parabolic partial differential equation.  
 Each node denoted by $A$ is a ``bottom'' node; \ie, all its edges are incoming. So $A$ might be an attractor.
 Each node denoted by $R$ is a ``top'' node; \ie, all its edges are out-going. 
 So $R$ might be a repellor. 
 Each node denoted by $S$ is a ``saddle'' node; \ie, it has in-coming and out-going edges.
 These graphs show how the dynamical system's nodes are dynamically related to each other, and for simplicity, they may ignore unbounded trajectories that come from or go to $\infty$ as $t\to\pm\infty$. 
 Or infinity may be one of the nodes.
\noindent
{\bf The Chafee-Infante (CI) PDE.} Its graph can have any finite number of levels depending on the parameter $\lambda$; see Eq.\eqref{eq: Chafee-Infante}.
}
 \label{fig:fiedler graphs}
\end{figure}

\black
\medskip\noindent{\bf $\eps$-chains for maps.}
Let $F$ be a map.
An $\eps$-chain from $x$ to $y$ under $F$ is a sequence of points $z_1,\dots,z_k$, \blue $k\geq2$, \black such that $z_1=x$, $z_k=y$ and
\beqn
\label{eq: chain condition}
d(F(z_i),z_{i+1})< \eps
\text{ for }i=1,\dots,k-1.
\eeqn

\medskip\noindent
{\bf $\eps$-chains for differential equations.}
Assume that $X$ is a Banach space with norm $\|\cdot\|$.
We provide an equivalent concept for differential equations that reduces to~\eqref{eq: chain condition} for the time-1 map of the flow.
Consider the differential equation
\beqn
\label{eq: x'}
\frac{dx}{dt}=g(x),
\eeqn
where we assume solutions $G^t(x)$, with $G^0(x)=x$, are uniquely defined for all forward time, $t\geq0$. 
We refer to $G$ as the flow of~\eqref{eq: x'}.

Let $z(t)$ be defined in $[0,T]$, where $T>0$.
We say that {\BF $z(t)$ is $\eps$-close to solutions $x(t)$} of \eqref{eq: x'} if, for every $t_0,t_1$ such that $0\leq t_0<t_1\leq T$, $t_1-t_0\leq1$,
\beqn
x(t_0)=z(t_0)\text{ implies }\|x(t_1)-z(t_1)\|<\eps. 
\eeqn
In practice, $z(t)$ will be continuous or piece-wise continuous.

Given $\eps>0$, we say that there is an {\BF $\eps$-chain from $x_0$ to $y_0$} if there is a curve $z(t)$ in $X$ defined on some interval $[0,T_\eps]$, where $T_\eps\geq1$, such that:
\begin{enumerate}
    \item $z(0)=x_0$ and $z(T_\eps)=y_0$;
    \item $z(t)$ is $\eps$-close to solutions of \eqref{eq: x'}.
\end{enumerate}
%
\blue
Our definition of $\eps$-chains for differential equations can actually be applied to all continuous-time dynamical systems and we leave the formulation to the reader. 
\black

\medskip\noindent
{\bf The downstream relation.} For a dynamical system, we say {\BF $y$ is downstream from $x$} (or {\BF $x$ is upstream from $y$}) if \blue either $y=x$ \black or there is a compact set $K$ such that, for every $\eps>0$, there is an $\eps$-chain from $x$ to $y$ lying in $K$.
This definition is for both discrete-time and continuous-time systems.
\blue
The downstream relation is transitive, that is if $z$ is downstream from $y$ and $y$ is downstream from $x$, then $z$ is downstream from $x$, and reflexive, that is every $x$ is downstream from itself.

\black
\medskip\noindent
{\bf Example of $\eps$-chains for differential equations.}
Assume that, for some $L>0$ and some bounded ball $B\subset X$, $g$ satisfies the Lipschitz condition 
\beqn
\label{eq: Lip}
\|g(x)-g(y)\|\leq L \|x-y\|\text{ for $x,y\in B$.}
\eeqn
Assume that, for some $T\geq1$ and $\eps>0$, $u(t)$ satisfies
\beqn
\label{eq: u(t)}
\int_t^{t+1} \|u(t)\|dt<\eps\cdot e^{-L}\text{ for each $t\in[0, T-1]$}.
\eeqn
Let $z:[0,T]\to B$ be a solution of 
\beqn
\label{eq: z'}
\frac{dz}{dt}=g(z(t))+u(t)
\eeqn
Let $x_0=z(0)$ and $y_0=z(T)$.
Then, by the Gronwall inequality, $z(t)$ is an $\eps$-chain from $x_0$ to $y_0$.

If $X$ has an inner product $\langle\cdot,\cdot\rangle$ then, instead of \eqref{eq: Lip}, we can require that
\beqn
\langle g(x)-g(y), x-y\rangle\leq L \langle x-y,x-y\rangle^{\nicefrac 1 2}\text{ for $x,y\in B$}.
\eeqn

\medskip\noindent
{\bf Chain-recurrent set.}
A set $C\subset X$ is chain-recurrent if, for every $x_0,y_0\in C$, $y_0$ is downstream from $x_0$ (and $x_0$ is downstream from $y_0$ and each $x_0\in C$ is downstream from $x_0$).

\medskip\noindent
{\bf Maximal chain-recurrent set.}
A chain-recurrent set is maximal if it does not lie in a strictly larger chain-recurrent set.
For instance, in Figure~\ref{fig:cr}, each limit set $N_1$, $N_2$ and $N_3$ is maximal. In contrast, for the logistic map $\ell_4(x)=4x(1-x)$ on $X=[0,1]$, it is well known that there are infinitely (countably) many periodic orbits, but none of them is maximal, since there are orbits of $\ell_4$ that are dense in the whole $[0,1]$ (recall that it is often the case that a chaotic attractor contains infinitely many periodic orbits).
In this case, the only maximal limit set of $\ell_4$ is $[0,1]$ itself.
As we show in Section~\ref{sec: proof}, maximal chain-recurrent sets are always closed and, when the system has a global trapping region, they are compact.

\medskip\noindent
{\bf Chain Graph.}
A graph of a dynamical system is a way to describe the relationships between some of its invariant sets.
The graph is abstracting information about the dynamical system.
It does not retain information such as the nature of nodes.
For instance, knowing that there is an attracting node does not tell us if the node is
a fixed point, a periodic orbit, a chaotic attractor, or something more exotic.
A researcher may then choose to add information as we have labeled the nodes in Figs.~\ref{fig:cr} and~\ref{fig:fiedler graphs}.

Graphs of dynamical systems are rarely discussed in the literature. 
Here, our concept of graph is based on a well-known concept, chains.
%

Any (directed) {\BF graph} consists o f a set of nodes and (directed) edges that run from one node to another. 
Hence, in order to define the chain graph of a dynamical system $F$, we only need to specify what its nodes and edges are.

\black
\medskip\noindent
{\bf Nodes of the chain graph.}
Each {\BF node of the chain graph} of $F$ is a maximal chain-recurrent set of $F$.

\medskip\noindent
{\bf Edges of the chain graph (chain edges)}. Let $A$ and $B$ be two distinct nodes of the chain graph of $F$. 
In our graph, there is at most one directed edge from $A$ to $B$.
We say that there is a (directed) {\BF edge from $A$ to $B$}, and write $A\Fto B$, if there are points $x\in A$ and $y\in B$ such that $y$ is downstream from $x$.
We say that the edge from $A$ to $B$ is a {\bf strong edge} if there is a two-sided trajectory  $\tau$ for which $\alpha(\tau)$ and $\omega(\tau)$ are non-empty and
\beq
\alpha(\tau) \subset A \text{ and }\omega(\tau) \subset B.
\eeq
We sometimes say an edge is {\bf weak} if it is not strong (see Conjecture~\ref{c2} in Section~\ref{sec: conj}).

\blue
\medskip\noindent
{\bf For the specialist:} the downstream relation is a pre-order. 
Its equivalence classes here are the maximal chain-recurrent sets. 
Thus, the chain graph, whose nodes are the equivalence classes and whose edges represent the pre-order, is a transitive acyclic graph. In other words, the graph represents the partial order relation between nodes.

\medskip\noindent
{\bf Chain graphs are topological spaces.} 
Convergence of sequences is a key concept in dynamical systems. 
In graph theory there is no natural concept of a sequence of nodes converging to a node but there is such a concept for chain graphs.
Indeed, we say that a sequence of nodes $N_k$ converges to a node $N_0$ if the minimum distance in $X$ between $N_k$ and $N_0$ goes to 0.

\medskip
From now on, {\bf the unmodified term \emph{graph} means \emph{chain graph}.}

\black
\medskip\noindent
{\bf Nodes and the global attractor.} In order for the graph to have at least one node (i.e. to be non-empty), there must be at least one non-empty limit set.
Hence, the graph and the global attractor are non-empty at the same time. 
Moreover, each node lies in the global attractor (see Proposition~\ref{p2} below) and, since the global attractor is closed, the closure of the union of all nodes also lies in the global attractor.

\medskip\noindent
{\bf Connected sets of chain-recurrent points and nodes.} In a metric space, if $Y$ is the union of nodes of a chain graph and $C$ is a connected set in $Y$, then $C$ is in a single node.
This is true because: 1. for any $\eps>0$ and any two points $x,y\in C$, there is an $\eps$-path from $x$ to $y$, \ie\ a finite sequence $p_1,p_2,\dots,p_n\in C$ such that $p_1=x$, $p_n=y$ and $d(p_i,p_{i+1})<\eps$ for all $i<n$; 2. given an $\eps$-path in $C$ from $x$ to $y$,
there is an $\eps$-chain from $x$ to $y$ which includes all the points of the $\eps$-path.

\medskip\noindent
{\bf Top and bottom nodes.}
\label{node}
We say that a node is a {\bf top node} if it has no incoming edges and is a {\bf bottom node} if it has no outgoing edges.
These are analogous to, but not equivalent to, the concepts of repellor and attractor.
A common occurrence would be a top node for which there is a trajectory coming from infinity to the top node. 
Hence, it would be not a repellor, and this unbounded trajectory is not reflected in the graph.

\medskip\noindent
{\bf Example of a graph.} In Figure~\ref{fig:cr}, every point $x$ in the region between $N_1$ and $N_2$ is such that $\omega(x)=N_2$ and $\alpha(x)=N_1$.
Similarly, every point $x$ in the region between $N_2$ and $N_3$ is such that $\omega(x)=N_2$ and $\alpha(x)=N_3$.
Hence, in the chain graph of this system, there are strong edges $N_1\Fto N_2$ and $N_3\Fto N_2$ and no other edges.

\medskip\noindent
{\bf Example of a graph.} The differential equation $x'=x(1-x^2)$ has three nodes, each of which is a fixed point: $N_1=\{-1\}$, $N_2={0}$, $N_3=\{1\}$.
However, in this case, the edges are reversed: $N_2\Fto N_1$ and $N_2\Fto N_3$.

\medskip\noindent
{\bf Example: a graph with weak edges.} An example of a graph with weak edges is the graph of the differential equation $x'=x(x-1)^2(x-2)^2$.
The nodes are 0, 1 and 2, and there are strong edges from 0 to 1 and from 1 to 2 and a weak edge from 0 to 2.

\medskip\noindent
{\bf Example: a graph without edges.} 
The chain graph of any dynamical system with a global trapping region has at least one node but might have no edges. This is the case, for instance, of the logistic map $\ell_4$ on $X=[0,1]$, whose chain graph has a single node $N=X$ and no edge.

\medskip\noindent
{\bf Example: the ``Connection Graph'' of Fiedler and Rocha.}
In~\cite{FR96,FR00,FR08,FR09,FR10,FRW12,FR21,FR22,FR23}, Fiedler and Rocha introduced a ``connection graph'' in their study of the global attractor of scalar semilinear parabolic partial differential equations (PDEs).
In Section~\ref{sec: PDE}, ``A PDE whose chain graph is completely described'', we describe their results concerning the
scalar reaction-diffusion PDE (where $x,u\in\R$),
\beqn
\label{eq: Chafee-Infante}
u_t = u_{xx} + \lambda u(1-u^2),\text{ where } \lambda>0.
\eeqn
Under their hypotheses, there is a global attractor, and each of its points is a constant function $u(t,x)$ that is independent of $x$.
Each node is a single-steady state point $u(t,x)$, which is independent of $t$ and $x$.
Each edge is a strong edge.
A few of their examples are illustrated in Fig.~\ref{fig:fiedler graphs}.
Their connection graph is a special case of our chain graph.

\blue
\medskip\noindent
{\bf Two important ``downstream'' properties.}
The next two propositions show the following two important facts.
First, independently on the topological properties of the phase space $X$, the qualitative dynamics of a dynamical system with a global trapping region coincides with the qualitative dynamics of a dynamical system on a compact space.

\medskip
\begin{proposition}
    \label{p2}
    Assume that a dynamical system $F$ has a global trapping region.
    Then no point outside of the global attractor is downstream a point in the global attractor.
    In particular, the chain-recurrent set of $F$ is contained in the global attractor and the graph of $F$ is the same as the graph of the restriction of $F$ to its global attractor.
\end{proposition}
\medskip
\begin{corollary}
    All maximal chain-recurrent sets of $F$ are compact.
\end{corollary}
\medskip

Moreover, the qualitative dynamics of a continuous-time system coincides with the one of its time-1 map.

\medskip
\begin{proposition}
    \label{p4}
    Assume that a dynamical system $F$ has a global trapping region.
    For any $T>0$, set $f=F^T$ its time-$T$ map.
    Then $f$ has the same chain-recurrent set, the same nodes, the same edges, and the same graph as $F$.
\end{proposition}
\medskip

Hence, the class of continuous-time systems with a global trapping region essentially coincide with the class of discrete-time systems with a compact global attractor.

\black
We postpone the proofs to Section~\ref{sec: proof}.

\medskip\noindent
{\bf Adjacent nodes.}
We call nodes $A$ and $B$ {\bf adjacent} if $A\Fto B$ and there is no node $C$ such that $A\Fto C$ and $C\Fto B$.
Assume $A, B$, and $C$ are nodes of a graph.
If there are edges
$A \Fto B$, $B \Fto C$ then there is an edge from $A \Fto C$.
Figures showing the graph will sometimes omit $A \Fto C$ to avoid having a cluttered figure and it will only show adjacent edges. 
For instance, that has been done systematically in Figure~\ref{fig:fiedler graphs}.

\medskip\noindent
{\bf Topology of the chain graph.}
In graph theory, the set of nodes of a graph is usually endowed with the discrete topology since, in general, there is no natural topology to give to this set.
Consider the graph where each point of $X=[0,1]$ is a node, and there are no edges. 
By graph definitions, it is not connected. 
For the chain graph of $X$, for every $\eps>0$, there is an $\eps$-chain between any pair of points.
Hence, $X$ has a single node consisting of the entire interval, so its graph is connected.

\medskip\noindent
{\bf Gradient dynamical systems.}
For dynamical systems that are ``gradient systems'' (see~\cite{Rob01}, Section~10.6.1), each node of the chain graph is either a fixed point, i.e. a constant steady state solution (\cite{Rob01}, Proposition~10.12), or a connected set of fixed points. 
When the set of fixed points is discrete, the global attractor is the union of the fixed points and some number of heteroclinic orbits (\cite{Rob01}, Proposition~10.13).
These heteroclinic orbits give rise to the edges of the corresponding chain graph.
These results have recently been generalized by
Lappicy for a wide class of quasilinear reaction-diffusion parabolic PDEs~\cite{Lap18} and fully nonlinear PDEs~\cite{Lap23}.

\medskip\noindent
{\bf Examples where $X$ is a circle.}
Consider now the simpler case of a continuous-time dynamical system on $X=\bS^1$ consisting of a fixed point (say the ``south pole'' $S$) and all other points moving clockwise.
Then, the global attractor $G$ is equal to the whole circle, and the chain graph consists of a single node $N=G$.
Consider now the case where there is a closed set $C$ of fixed points in the circle, and every other point moves clockwise.
Even in this case, $G=\bS^1$, and there is a single node $N=G$. 

\medskip\noindent
{\bf Connectedness of the chain graph.}
\blue In traditional graph theory, \black a graph is connected if there is a (not necessarily oriented) path from any node to any other node in the graph.
While ``paths'' are allowed to have infinitely many edges, \blue in absence of a topology, \black traditional graph theory \blue is not able to tell \black what is meant by saying a path with infinitely many edges connects node 1 to node 2. 
\blue 
The presence of a topology on chain graphs makes it possible to give a natural meaning to this concept.

More specifically, we define connectedness for a chain graph as follows.
We say that a collection of nodes of a chain graph is {\bf closed} if it is a closed set with respect to the topology on the graph.
Then, we say that the {\BF\blue chain graph $\Gamma$ is connected} when, for each decomposition of its nodes into two closed disjoint collections $\cM$ and $\cN$, there is an edge from some node of $\cM$ to some node of $\cN$ or vice versa.

Equivalently, denote by $M\subset X$ and $N\subset X$ the corresponding sets of node points of $\cM$ and $\cN$ respectively.
Then $\Gamma$ is connected when, for each decomposition of its nodes into two disjoint collections $\cM$ and $\cN$ such that $M$ and $N$ are closed and disjoint, there is an edge from some node of $\cM$ to some node of $\cN$ or vice versa.

\begin{figure}
    \centering
\tikzset{
        set arrow inside/.code={\pgfqkeys{/tikz/arrow inside}{#1}},
        set arrow inside={end/.initial=>, opt/.initial=},
        /pgf/decoration/Mark/.style={
            mark/.expanded=at position #1 with
            {
            \noexpand\arrow[\pgfkeysvalueof{/tikz/arrow inside/opt}]{\pgfkeysvalueof{/tikz/arrow inside/end}}
            }
        },
        arrow inside/.style 2 args={
            set arrow inside={#1},
            postaction={
                decorate,decoration={
                    markings,Mark/.list={#2}
                }   
            }
        },
        -<-/.style={decoration={markings,mark=at position 0.5 with {\arrow[red,>=stealth,scale=1.5]{<}}}, postaction={decorate}},
        ->-/.style={decoration={markings,mark=at position 0.5 with {\arrow[red,>=stealth,scale=1.5]{>}}}, postaction={decorate}},
        -<--/.style={decoration={markings,mark=at position 0.5 with {\arrow[red,>=stealth]{>}}}, postaction={decorate}},
    }

\definecolor{peg}{rgb}{1, 0.8, 0.1}
\setlength{\tabcolsep}{10pt} 
\renewcommand{\arraystretch}{1.5} 

  \begin{tikzpicture}[every node/.style={scale=1.5}]
    \pgfmathsetmacro{\goldenRatio}{(1+sqrt(5)) / 2}

    \draw [color=black,->-,thick] (7.5,0) -- (10,0);
    \draw [color=black,-<-,thick] (5.5,0) -- (7.5,0);
    \draw [color=black,->-,thick] (4,0) -- (5.5,0);
    \draw [color=black,-<-,thick] (2.75,0) -- (4,0);
    \draw [color=black,->-,thick] (1.75,0) -- (2.75,0);
    \draw [color=black,-<-,thick] (1.,0) -- (1.75,0);
    \filldraw [black] (0,0) circle (2pt);
    \filldraw [black] (1,0) circle (2pt);
    \filldraw [black] (1.75,0) circle (2pt);
    \filldraw [black] (2.75,0) circle (2pt);
    \filldraw [black] (4,0) circle (2pt);
    \filldraw [black] (5.5,0) circle (2pt);
    \filldraw [black] (7.5,0) circle (2pt);
    \filldraw [black] (10,0) circle (2pt);
    \node at (0.5,0) {$\dots$};
    \node at (0,-0.5) {$0$};
    \node at (7.5,-0.5) {$\frac{1}{2}$};
    \node at (5.5,-0.5) {$\frac{1}{3}$};
    \node at (4,-0.5) {$\frac{1}{4}$};
    \node at (2.75,-0.5) {$\frac{1}{5}$};
    \node at (1.75,-0.5) {$\frac{1}{6}$};
    \node at (1,-0.5) {$\frac{1}{7}$};
    \node at (10,-0.5) {$1$};

  
  \end{tikzpicture}    
    \caption{A dynamical system on $X=[0,1]$ with a connected chain graph where a node (the fixed point at 0) has no ingoing nor outgoing edges.
    Points 0 and $1/k$, $k=1,2,\dots$, are fixed while points between $1/k$ and $1/(k+1)$ move in the direction indicated by the red arrows.
    The picture is not to scale.
    The graph of this system essentially coincides with the picture above, where the segment between $1/k$ and $1/(k+1)$ represents now an edge between those two nodes. 
    The direction of the edge is the one indicated by the red arrow.
    }
    \label{fig: connected graph}
\end{figure}

\black
\medskip\noindent
{\bf An example with countably many nodes.}
A simple case where the definitions of graph connectedness disagree is given by the flow of the ODE $\dot x=x\sin\frac{\pi}{x}$ for $x\in[0,1]$ (see Figure~\ref{fig: connected graph}).
\blue 
Here, $X$ is the interval $=[0,1]$
\black
and $x$ is a node if it is a stationary point, namely $x=1/n$, $n=1,2,\dots$, or $x=0$.
The node $1/n$ is an attractor if $n$ is odd and a repellor if $n$ is even.
There is a directed edge from $1/n$ to $1/m$ if and only if $n$ is even and $m=n\pm1$.
\blue
There are no other edges.
In the definition of the connectedness of a chain graph, we are asked to partition the nodes into two disjoint closed subsets.
If that is done for this example, then one of the sets cannot be the single node $0$ since remaining nodes would not constitute a closed set. 
Hence, the chain graph of this ODE is connected.
On the other side, node 0 has no edges.
Hence, by the traditional definition, this graph is not connected.

\black
\smallskip\noindent
{\bf Another example with countably many nodes: the Feigenbaum parameter value.}
Feigenbaum discussed dynamical systems depending on a parameter $r$, in which there is a sequence of bifurcations at parameter values $r_n$, where $r_n$ converges to some point $r_\infty$.
As in the logistic map, it is frequently true that, at $r_\infty$, there is a Cantor set $C_\infty$ that attracts almost all initial conditions near $C_\infty$, and there is a sequence of unstable periodic orbits $P_n$ of period $2^n$ that approach $C_\infty$ as $n\to\infty$.
From the point of view of chain graphs, $C_\infty$ is a node as is each $P_n$.
Note that the node $P_n$ consists of the $2^n$ points of the periodic orbit.
For the logistic map, there is an edge from $P_{n+1}$ to $P_n$, yielding the following picture:
\beq
P_0\mapsto P_1\mapsto\cdots\mapsto P_n\mapsto P_{n+1}\mapsto\cdots C_\infty
\eeq
The above picture shows only the edges between adjacent nodes, and there may be many other nodes in addition. 
This set of edges is a path from $P_0$ to $C_\infty$.
While this path is a directed path, with each edge leading toward $C_\infty$, in general the definition of connectedness allows for a path of edges without a consistent direction.




\smallskip
We are now ready to state our main result.

\black
\vskip 0.3cm
\noindent
{\bf Theorem~\ref{Connectedness of graph}}({\bf Connectedness of chain graphs}) 
    Assume the global attractor of $F$ is connected. 
    Then the chain graph of $F$ is connected.

\medskip\blue
The assumption that the global attractor is connected is almost always satisfied in practice. 
For example, by Proposition~\ref{p1}, if $F$ has a connected trapping region, then the global attractor is connected.
In practice, trapping regions are compact balls and, in infinite dimensions such as for certain dissipative PDEs, the compact ball is in some subspace.
See~\cite{GobS97} for some results on the connectedness of global attractors under stronger hypotheses.

\medskip\black
Sections~\ref{sec: logistic map}-\ref{sec: PDE} provide four examples with details.
This theorem applies to the wide class of dynamical systems that have a connected global attractor, including all systems that have a connected trapping region. 
Notice that our theorems do not assume that $X$ is complete or connected.
We postpone to Section~\ref{sec: proof} the proof of the theorem. 

Here is a non-trivial application of our theorem. 
Recall that Sheldon Newhouse \cite{New74} showed that, for diffeomorphisms of the plane, it is quite common for a map that is chaotic for one parameter value to have uncountable many nearby parameter values for which the map has infinitely many attractors coexisting within some compact set. 
Newhouse's results, though, do not say anything about non-attracting nodes. 
Since attracting nodes have no outgoing edges, our Theorem~\ref{Connectedness of graph} implies that Newhouse systems must have at least one non-attracting node.
Many planar examples satisfy Newhouse's criterion that there is a tangency between stable and unstable manifolds. 

One might think that, if there are infinitely many periodic attractors in a compact set, the basins must be small. But the basins are not restricted to that compact set.
It is not uncommon for examples of maps $F$ to have a Jacobian determinant $c\bydef|\det DF(x)|<1$ that is constant in phase space, for example, the Henon map or the forced damped pendulum. 
If $B$ is a basin and its area is $\beta$, then $F(B)=B$ and area$(B)=$area$(F(B))=c\cdot$area$(B)$.
Hence, $$\beta=c\beta,$$ namely $\beta=\infty$ since $\beta>0$.
Hence, in such cases, even when there are infinitely many basins, each one has infinite area.

\medskip\noindent
{\bf Morse graph.} 
Morse described the graphs of certain special ``gradient flows''. 
Generalized Morse graphs are discussed by Smale~\cite{Sma67}, Norton~\cite{Nor95}, Osipenko~\cite{Osi06}, Mizin~\cite{Miz02}, Yokoyama~\cite{Yok20}, Mischaikow {\em et al.}~\cite{VM22} using a variety of formulations, often inspired by Morse,  Conley's $\eps$-chains and Auslander's prolongations.
There are very few other papers that define what a graph of a general dynamical system is~\cite{Nor95,FR96,Miz02,Osi06,DLY20,DLY21}. 
Smale, for example, only defined a graph in the case of ``Axiom-A'' diffeomorphisms. 
In that case, each node is an invariant set with a dense trajectory, and there is an edge from node $N$ to node $M$ if the unstable manifold of $N$ intersects the stable manifold of $M$.
\begin{figure}
    \centering
    \includegraphics[width=10cm]{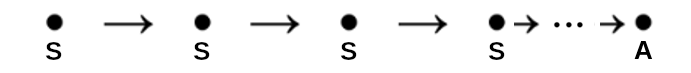}
    \caption{{\bf The Morse graph of the logistic map.} 
    For the logistic map $f(x)=ax(1-x)$, for $r\in(1,4]$ and $x\in[0,1]$,
    The Morse graph of each logistic map is a ``tower''~\cite{DLY21}. Namely, there is an edge between each pair of nodes. 
    Depending on the parameter, the tower can have any finite number of nodes or an infinite number that converge to the unique attracting node, denoted by $A$. 
    When the number of nodes is infinite, the attractor is quasiperiodic.
    Otherwise, the attractor is either periodic or chaotic.
    Each saddle node, denoted by $S$, is either a periodic orbit or a Cantor set saddle with chaotic dynamics. 
    {\bf Reduced graphs:} In our figures, if there is an edge from a node $N_1$ to a node $N_2$ and one from $N_2$ to node $N_3$, then we do not draw the edge from $N_1$ to $N_3$.
}
    \label{fig: logistic map}
\end{figure}
When Smale's graph is defined, it is a chain graph. 
Different definitions of general recurrence can yield different Morse nodes and may have additional edges.

\black

\section{An example in dimension 1: the logistic map}
\label{sec: logistic map}
We have good news and bad news about the logistic map $\ell_r(x)=r x(1-x)$, $X=[0,1]$, and its graph.
For $r\in(1,4]$, $X$ is forward invariant, and so it has a graph. 
The good news is that in~\cite{DLY20},  its chain graph has been determined for every $r$.
There are uncountably many $a$ values for which the graph has infinitely many nodes, and for every $r$, there is a directed edge between each pair of nodes. 
For each pair of nodes $M$ and $N$, there is either an edge from $M$ to $N$ or $N$ to $M$ but not both.
As already stated, there are no loops in the Morse graph.
We call such a graph a {\bf tower} (see the rightmost figure (T) in Fig.~\ref{fig:fiedler graphs}).

Usually, the tower has finitely many nodes.
A technical point of our proof depends upon the shortest distance $d_N$ between a node $N$ and the critical point $c=0.5$.
We show that if $d_M>d_N$, then there must be a strong edge from $M$ to $N$.

The simplest example of an $r$ for which there are infinitely many nodes is the Feigenbaum value, as we discussed above, $r_\infty=3.5699...$.
At that value, there are infinitely many nodes $M_n$, where $M_n$ is the node consisting of an unstable periodic orbit of period $2^n$.
There is also the attractor node $M_\infty$.
Our logistic theorem says that, in this case, there is an edge from $M_i$ to $M_j$ if and only if $i<j$, where $j$ can be $\infty$.

\black
\begin{figure}
    \centering
    \includegraphics[width=\linewidth]{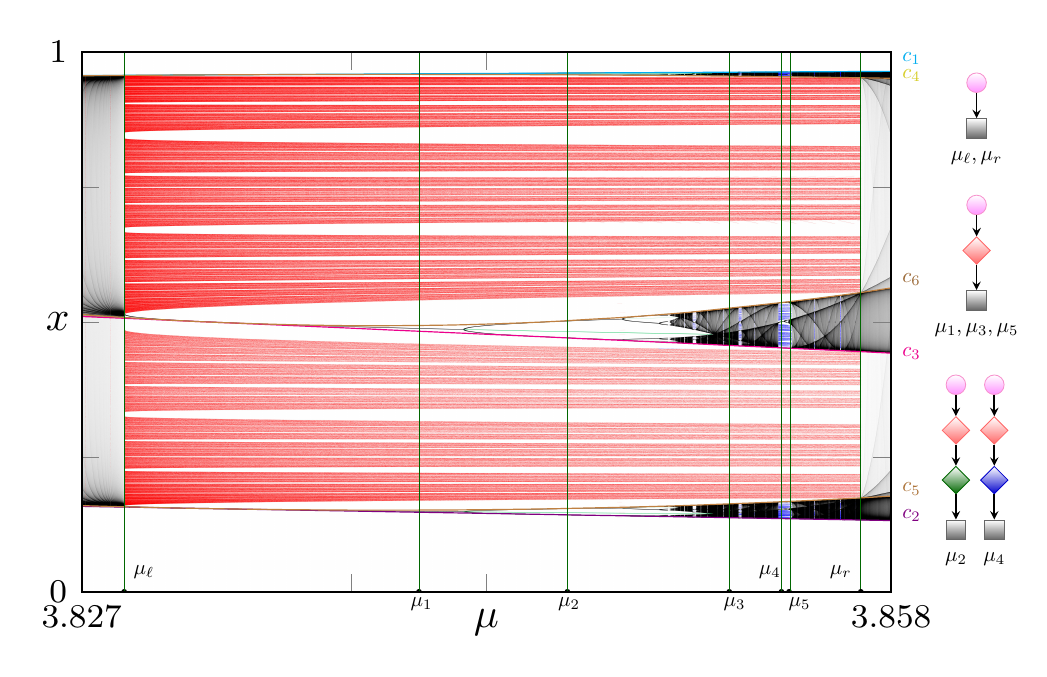}
    \caption{{\bf Bifurcation diagram and chain graph of the logistic map within the period-3 window.} This bifurcation diagram shows some of the nodes of the logistic map in the period-3 window, from~\cite{DY25}.
    Parameters $\mu_\ell$ and $\mu_r$ denote, respectively, the parameter values of the left and right endpoints of the window. At the right of the picture are shown the chain graphs of the map at the parameter values marked in the horizontal axis. Circle nodes represent repelling nodes (in this case, the $x=0$ fixed point). Diamond nodes represent saddles, colors correspond to the color of the saddle in the bifurcation diagram -- red and blue saddles are chaotic Cantor sets, the green saddles are periodic orbits. Square nodes represent attracting nodes.}
    \label{fig: lm}
\end{figure}

A caution: discussions of Feigenbaum-like results emphasize the bifurcation diagram as $r\to r_\infty$ from the left, and the map undergoes period-doublings as the attracting periodic orbit's changes from period $2^n$ to period $2^{n+1}$.
Our result is not about this bifurcation sequence. 
It is about the parameter $r_\infty$ and the periodic orbits at $r_\infty$.
The attractor node $M_\infty$ is a complicated quasiperiodic orbit.
Perhaps we should call $r_\infty$ the first Feigenbaum parameter value, since there is a countable infinity of $r$ values corresponding to limits of period-doubling bifurcation sequences. 

More good news: logistic map graphs can be described completely, and, for some parameters, the graph has infinitely many nodes.
Choose any infinite sequence of the labels ``chaotic saddle'' node and ``periodic saddle'' node for the nodes $M_1,M_2,M_3,\dots$, where $M_1$ is the top node of (T) in Fig.~\ref{fig:fiedler graphs}, $M_2$ is the next node, etc.
Our theorem says that for each such sequence, there is an $r$ that has such a graph.
Since there are uncountably many ways to choose the sequence of chaotic and periodic nodes, it follows that the logistic map has uncountably many $r$ values having a graph with infinitely many nodes.

To the best of our knowledge, our logistic paper is the first to describe a graph with infinitely many nodes for a known dynamical system. 

The bad news: our ability to prove that the logistic map can have graphs with infinitely many nodes is not because the logistic map is unusually complicated, but rather we can prove our results because it is especially simple.
We expect every dynamical system that has cascades of period-doublings will have parameter values for which the graph has infinitely many nodes with a subgraph that is a tower, but such dynamical systems might have many other nodes, including many attractors. 
They would be more complicated than the logistic map graph, which has only one attractor for each $r$.

\black
\section{An example in dimension 2: the forced damped pendulum}
\label{sec: pendulum}
As we have argued in this article, the graph of even the simplest discrete-time dynamical systems on the line can have infinitely many nodes in their graph while having only one attractor node.
By the Newhouse theorem, the situation can be extremely more complicated in dimension 2 or more since there can be infinitely many attractor nodes.
The forced damped pendulum is such a case, but here we present a case where there are only three nodes.
Tiny changes in parameters can result in infinitely many attractors, but we do not know which perturbation would give that (see~\cite{SY15} and~\cite{SY15web}).

\begin{figure}
    \includegraphics[width=0.9\textwidth]{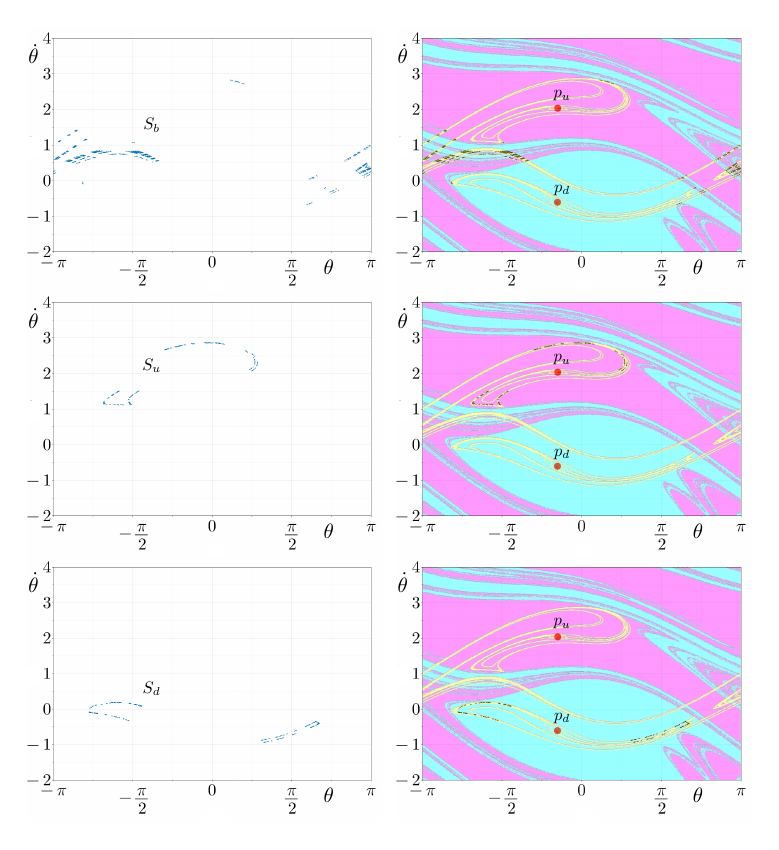}
    \caption
    {{\BF The three chaotic saddles of the forced damped pendulum.}
    Each panel shows the phase-space for the time-2$\pi$ Poincaré return map for parameters $\gamma = 0.2, \rho = 2$ (see (\ref{eqn:FDP})). 
    The global attractor is shown in yellow.
    The left column shows three chaotic saddles, and the corresponding panel on the right shows them embedded in a picture of the two basins (magenta and cyan) of attraction of the two fixed points (red dots) of the Poincar\'e return map. 
    Of course, chaotic saddles and fixed points lie in the global attractor.
    }
    \label{fig:DDP_saddles}
\end{figure}

\medskip\noindent
{\bf Numerical basin.} For a dynamical system, an attractor $A$ has a basin $B$ equal to the set of $x$ that are attracted to $A$.
The set $B$ is open when $A$ attracts a neighborhood of $A$.
Numerical determinations of the basin are often based on the examination of a grid of points.
There can be sets of measure zero, either periodic or chaotic, and their stable manifolds that are not in $B$ but are in the interior of the closure of $B$. Since they have measure zero, any finite grid is unlikely to obtain any points in such a set.
Hence, we define the {\bf numerical basin} of $A$ to be the interior of the closure of $B$.
As we have just described, the numerical basins can contain slightly more than the true basin. 
The numerical basin can include sets of measure zero that are not attracted by $A$, and so are not in the (true) basin of attraction.
Such sets are found in the following example.

We consider here the ODE
\begin{equation}
                \begin{split}
                    \ddot\theta  + \gamma \dot\theta 
                    + \sin\theta
                    = \rho \cos t,
                \end{split}
                \label{eqn:FDP}
\end{equation}
where $\gamma = 0.2$ and $\rho=2$.
We consider $\theta\in[-\pi,\pi] \bmod 2\pi$, and the derivative $\dot\theta\in\R$.
This ODE is periodically forced with period $t=2\pi$.
Let $f(\theta,\dot\theta)$ be the time-$2\pi$ map of the flow.
It can also be referred to as the Poincar\'e return map for the cross section at $t=0$ ($\bmod 2\pi$).
It has the same graph as the continuous flow on the 3-dimensional $(t,\theta,\dot\theta)$-space.
Then the map $f$ maps the cylinder $[-\pi,\pi]\times\R$ to itself. 
Since $\gamma>0$, $f$ is dissipative and, if $Y$ is sufficiently large, the compact cylinder $[-\pi,\pi]\times[-Y, Y]$ is mapped into itself, that is, it is an absorbing region and therefore has a compact connected global attractor.
If a bounded set $Q$ has area $q$, then its image $f(Q)$ has area $e^{-2\pi\gamma}q\sim 0.28q$, since the Jacobian of the flow in the $\theta$-$\dot\theta$ space has trace $-\gamma$.
Hence the time-$2\pi$ map shrinks areas by a factor of $e^{-2\pi\gamma}$.

All of our results for this map are numerical.
We cannot rule out the presence of additional tiny basins of attractions.
Hence, our graph represents the relationships between the nodes we have found.

The map $f$ has the following two attracting fixed points:
\beq
p_u \sim (-0.472, 2.037),\ p_d \sim (-0.478,-0.608).
\eeq

\begin{figure}
    \centering
    \includegraphics[width=4cm]{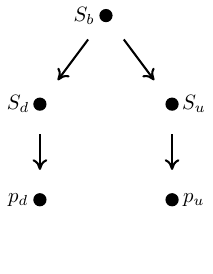}
    \caption{{\BF The graph of the forced damped pendulum.}
    The forced damped pendulum consists of five nodes, the fixed points $p_u,p_d$ and the three chaotic saddles $S_b$,$S_d$,$S_u$ shown in Fig.~\ref{fig:DDP_saddles}.
    Saddle $S_b$ is a top node and has edges to each of the other four nodes. 
    Saddle $S_u$ has a single edge to $p_u$, and similarly for $S_d$.
    Only adjacent edges of the graph are shown in the picture above.
    }
    \label{fig: DDP graph}
\end{figure}
In Fig.~\ref{fig:DDP_saddles}, the numerical basins of attraction of these two fixed points are plotted in cyan and pink. 
Then, the global attractor is added in yellow (recall that every node is a subset of the global attractor).
Next, the two fixed points $p_u$ and $p_d$ are plotted as large red dots.
Each of the fixed points is a node. 
In addition, there are
three more saddle nodes, painted black, each of which is
an invariant Cantor set (see~\cite{BGOYY88,NY89, NY91}).
We labeled these saddle nodes by $S_b$ (for boundary), $S_u$ (for upper) and $S_d$ (for down).
The first, $S_b$, lies in the boundary between the two basins, and the closure of its unstable manifold coincides with the global attractor.
Hence, there is an edge from $S_b$ to every other node of the graph (see Fig.~\ref{fig: DDP graph}).
 Each of the other two, $S_u$ and $S_d$, lies within a numerical basin, $S_u$ in the basin of $p_u$ and $S_d$ in the basin of $p_d$.
Since they are in a numerical basin, an arbitrarily small perturbation of $S_u$ gives a point in the basin of $p_u$, so there is an edge from $S_u$ to $p_u$, and similarly for $S_d$ and $p_d$.

Each of the three Cantor sets mentioned above 
is one of the following two types: 
\begin{enumerate}
    \item \textbf{Boundary Chaotic Saddle.} 
    This node is a ``top'' node; \ie, all its edges are out-going.
    This set is located on the basin boundary and is chaotic. 
    It is also an invariant limit set, with unstable directions pointing towards the basins. 
    For a discussion on the dynamics of this set, see~\cite{BGOYY88}. 
    \item \textbf{Embedded Chaotic Saddle.} 
    There are two distinct invariant chaotic Cantor sets in the interior of the closure of the two basins of attraction. 
    Recall that these saddles have previously been reported by H. Nusse and J. Yorke in papers for the PIM triple procedure~\cite{NY89, NY91}. 
\end{enumerate}
\section{An example in dimension 3: the Lorenz system}
\label{sec:Lorenz}

\begin{figure}
    \centering\includegraphics[width=14cm]
    {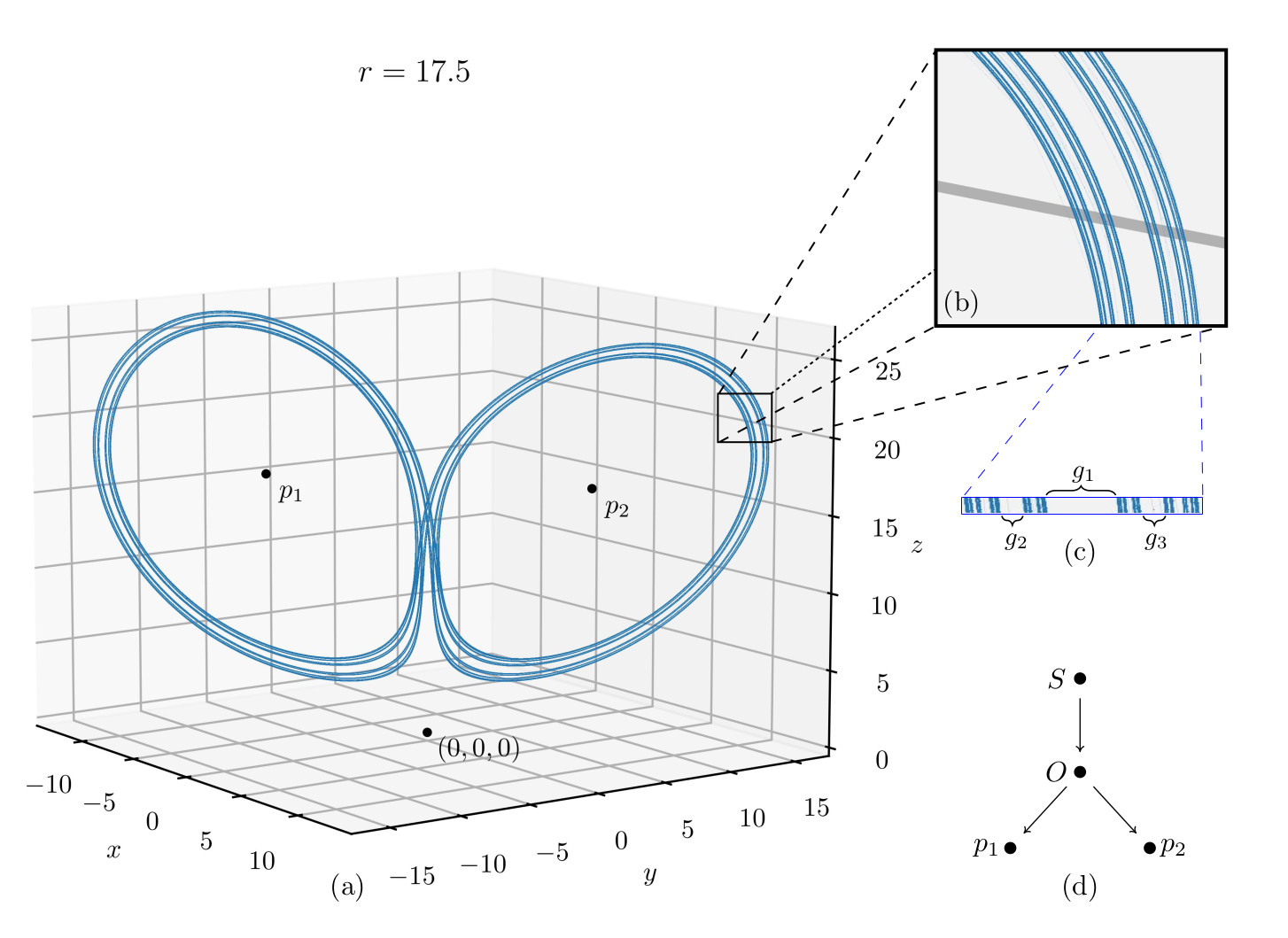}
    \caption{{\bf All nodes of the Lorenz system for $r=17.5$.} 
    (a): The Lorenz chaotic saddle is a node (CS); (b): blowup of a portion of CS; 
    (c): cross-section of the blowup that looks like a middle-third Cantor set.
    On each side of the gap $g_1$, there are gaps $g_2,g_3$ about a third the size of $g_1$. 
    (d): graph of the system, where $O$ is the origin and $p_1,p_2$ are the attracting fixed points. 
    }
    \label{fig: Lorenz}
\end{figure}

\begin{figure}
\centering\includegraphics[width=10cm]
    {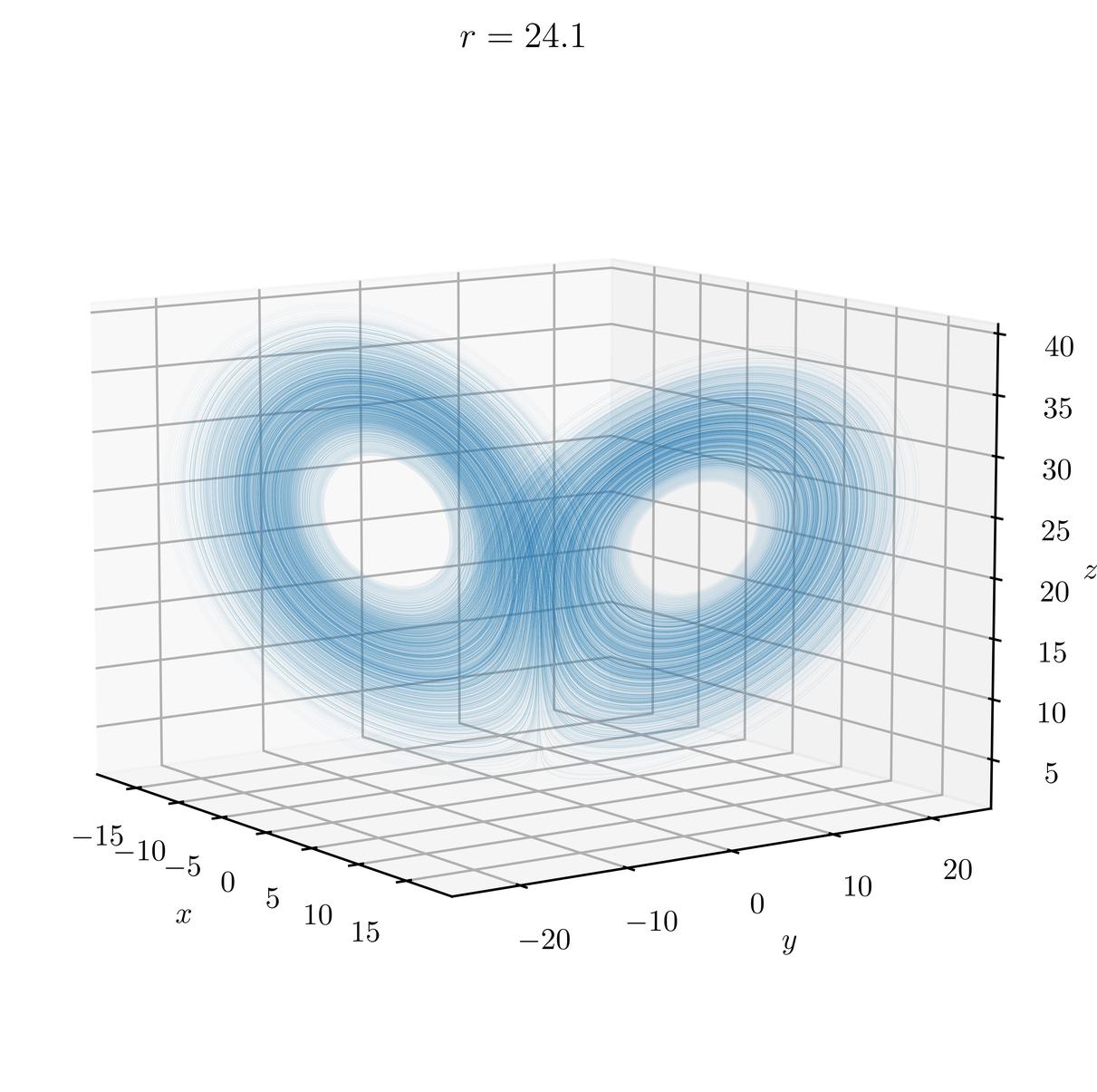}
    \caption{{\bf The chaotic attractor at $r=24.1$.}}
    \label{fig:lor24}
\end{figure}

In his seminal article~\cite{Lor63}, Lorenz introduced his famous three equations 
\begin{equation}
                \begin{split}
                    &x' = -\sigma x + \sigma y\\
                    &y' = -x z + r x - y\\
                    &z' = \phantom{-}x y - b z,
                \end{split} 
                \label{eqn:lorenz}
            \end{equation}
and investigated the butterfly-shaped attractor that occurs for $r=28.$ Throughout this section, $ \sigma = 10,\, b = 8/3$. Figure \ref{fig:lor24} shows the attractor for $r=24.1$.

Using a Lyapunov-like function, he showed there are balls in $\R^3$ that are forward invariant and eventually absorb every trajectory.
These balls are trapping regions. 
Hence, there is a global attractor. 

Later, it was discovered that for the same values of $\sigma$ and $b$, and $13.92\leq r \leq 24.06$, the chaotic attractor is replaced by a ``metastable'' set, which is now usually referred to as a chaotic saddle~\cite{CS82,ABS77,KY79,DLY21,YY79}. See Fig.~\ref{fig: Lorenz} (a)-(c) for two views of the chaotic saddle, S, at $r=17.5$. S is a saddle node in the graph of the system. 

In addition to the chaotic saddle, there are two fixed point attractors (sinks) in the system at $r=17.5$, labeled as $p_1$ and $p_2$, each of which corresponds to a distinct node in the graph. The coordinates of $p_1$ and $p_2$ are given by 
\beq
p_1 = \left(-\sqrt{b(r-1)},-\sqrt{b(r-1)}, r-1\right) \approx (-6.633, -6.633, 16.5)
\eeq
\beq
p_2 = \left(\sqrt{b(r-1)},\sqrt{b(r-1)}, r-1\right) \approx (6.633, 6.633, 16.5)
\eeq

Finally, the origin $O$ is a fixed-point saddle and hence a saddle node in the graph.

Regarding the structure of the graph, the stable manifold of the origin $O$ is contained in the boundary between the basins of attraction of $p_1$ and $p_2$~\cite{MGOY85}. So, the origin $O$ lies on the basin boundary and hence there is an edge from $O$ to both $p_1$ and $p_2$. It can also be noted that the basin boundary is a fractal. It occurs in the form of a Cantor set of planes. 

\begin{figure}
    \centering
    \includegraphics[width=0.75\linewidth]{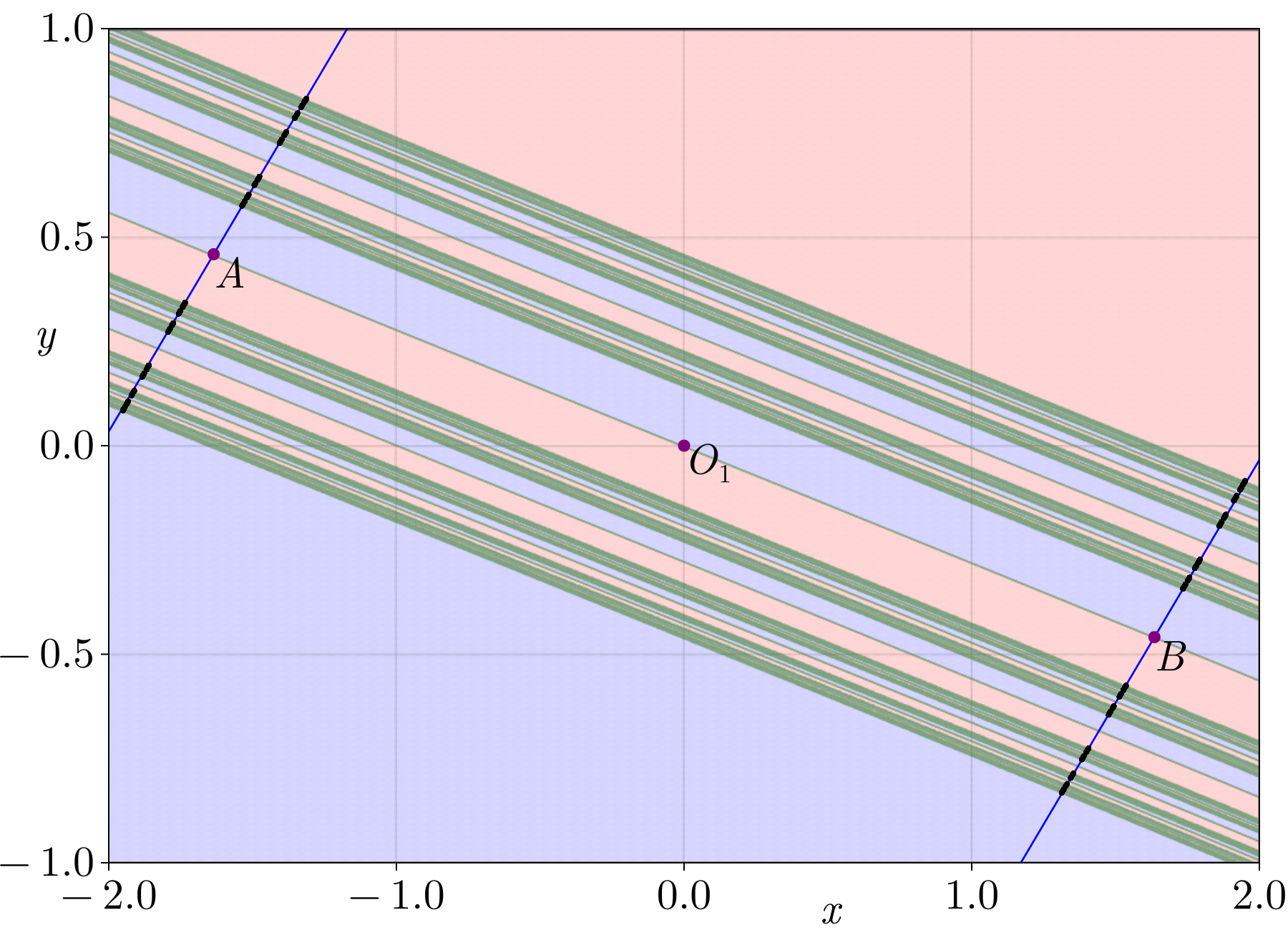}
    \caption{{\BF Poincaré return map image of the chaotic saddle at $r=17.5$, with Poincaré section $z=16.5$.} The point $O_1$ is $(0,0,16.5)$, which is trivially on the stable manifold of the origin.
    The isolated green lines, such as the line passing through $O_1$ are also on the stable manifold of the origin.
    The blue region is the numerical basin of $p_2$, and the red region is the numerical basin of $p_1$. All green lines represent the boundary between the basins. There is a Cantor set of green lines, each of which is a limit of the isolated green lines, and they are on the stable manifold of the chaotic saddle.   The cross section of the Cantor set of the chaotic saddle $S$ is plotted in black. 
    The dark lines through the cross section of $S$ are the cross section of the unstable manifold of $S$. The unstable manifold intersects the stable manifold of the origin at heteroclinic points, like $A$ and $B$. }
    \label{fig:poincarelorenz}
\end{figure}

Next, look at the Poincaré section $z=16.5$, containing both $p_1$ and $p_2$ in fig. \ref{fig:poincarelorenz}. The plot displays the intersection of the numerical basins of $p_1$ (in blue) and $p_2$ (in red) with the plane $z=16.5$. 
The cross section of the boundary between the basins, which appears as a Cantor set of lines, is plotted in green. 
The cross section of $S$, a different Cantor set, is plotted in black. 
Note that $S$ is contained inside the boundary between the basins, which implies there is an edge from $S$ to $p_1$ and $p_2$. 
Moreover, the unstable manifold (drawn in dark blue) of $S$ crosses the stable manifold of $O$, so that
there is an edge from $S$ to $O$. 
As mentioned in the caption of Fig. \ref{fig:poincarelorenz}, the lines passing through the Cantor sets are the cross sections of the unstable manifold of $S$.
These lines correspond to the stretching direction of the horseshoe analysis shown in~\cite{KY79}. 
The discussion above justifies the structure of the graph drawn in Fig. \ref{fig: Lorenz}(d).

\begin{figure}
    \centering
    \includegraphics[width=0.75\linewidth]{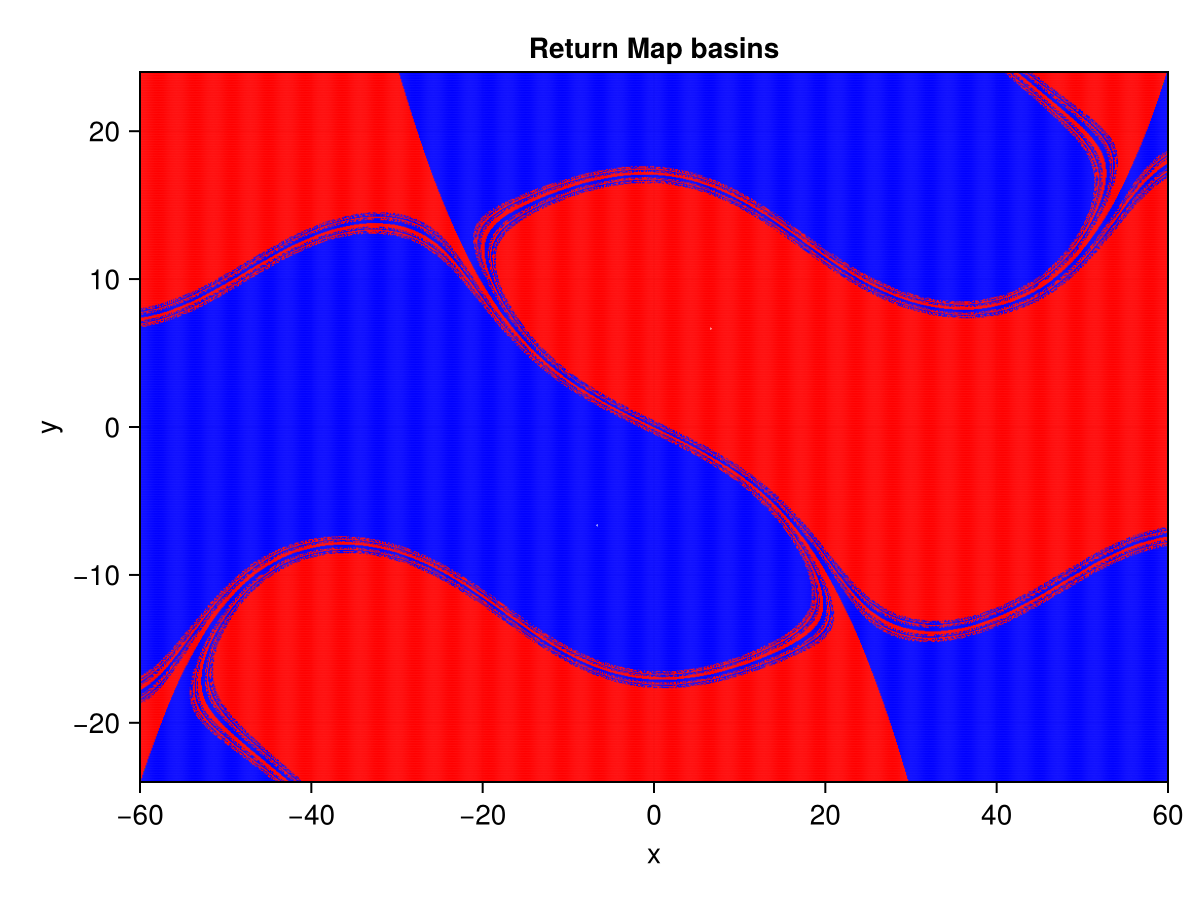}
    \caption{{\BF Basins of attraction of $p_1$ (red) and $p_2$ (blue) under the Poincaré return map with Poincar\'e section for $r=17.5$ and $z=16.5$}. The basin boundary is not a straight line, as it might seem in Fig. \ref{fig:poincarelorenz}. The image is also consistent with the numerical image of the stable manifold of the origin obtained by Doedel, Et. Al \cite{DKO15} for $r=18$, which has the same global dynamics as our test case of $r=17.5$.}
    \label{fig:basinslorenz}
\end{figure}

We chose the value $r=17.5$ so that the Cantor set looks like a standard middle-third Cantor set.
The largest gap's width-fraction is approximately 30\% here, and the fraction goes to 0 as $r$ increases towards 24.06. 
For $r>24.06$, the saddle becomes a chaotic attractor. 
A large-scale portrait of the basins is plotted in Fig. \ref{fig:basinslorenz}. 
This image is consistent with the numerically obtained images of the stable manifold of the origin in~\cite{Osi18} for $r=18.0$, which share the same qualitative global dynamics of our case.
We also used the software developed by Doedel, Krauskopf and Osinga in~\cite{DKO15} to generate the stable manifold of the Lorenz system at $r=17.5$ at the origin and verify the accuracy of our Figures~\ref{fig:poincarelorenz} and~\ref{fig:basinslorenz}.
Its results are in perfect agreement with our own.

Our numerical investigations in~\cite{DLY21} show that, between $r\sim$ 24.06 and $r\sim$ 30.1, the graph has three nodes (an attractor and two unstable fixed points) and is identical to the graph shown on the right side of Figure~\ref{fig:cr}, and that there are values of $r$ beyond 30.1 for which the graph has infinitely many nodes.
That article gives strong numerical evidence that the 3-dimensional Lorenz system has a behavior similar to that of the logistic map when Lorenz's parameter $r$ is near 208.
\black
\section{An example in infinite dimension: the Chafee-Infante PDE}
\label{sec: PDE}
Finding the complete graph of an infinite-dimensional system is, in general, a hopeless task.
Because of this, infinite-dimensional dynamical systems whose complete graph can be found are particularly precious.
Here, we consider one such system, and we use it to illustrate how, more generally, one can write a reaction-diffusion PDE whose graph is identical to the graph of a given ODE flow in finite dimension.

\medskip\noindent{\bf The ODE.}
Consider the ODE
\beqn
\label{eq: ODE}
\dot u = u( 1 - u^2)
\eeqn
in $X=\R$, whose global attractor is $[-1,1]$.
This ODE has three fixed points, namely two attractors $u=-1,1$ and a repellor $u=0$, and it has two (strong) edges. 
One is the heteroclinic orbit from 0 to 1, and the other is the heteroclinic orbit from 0 to -1.

\medskip\noindent
{\bf The PDE.}
Now
we consider the scalar reaction-diffusion PDE (where $x,u\in\R$),
\beqn
\label{eq: Chafee-Infante}
u_t = u_{xx} + \lambda u(1-u^2),\text{ where } 0<\lambda<1.
\eeqn
This is known as the Chafee-Infante PDE~\cite{CI74}.
We set periodic conditions on Equation~(\ref{eq: Chafee-Infante}), namely, $u(t,x)$ is periodic in $x$ with period $2\pi$. 
To our knowledge, the global attractor of this particular case has not been discussed in the literature.
We denote the circle by $\mathbb S^1$.
There are solutions of~(\ref{eq: Chafee-Infante}) for which $u_{xx}\equiv0$ and, for these solutions, Equation~(\ref{eq: Chafee-Infante}) becomes Equation~(\ref{eq: ODE}); see Theorem~\ref{thm:CI}.

The results of James Robinson in Chapter 11 of~\cite{Rob01} show that, under a general condition satisfied by $f(u)=u(1-u^2)$, the PDE above defines a semi-flow $F$ on $X=L^2(\mathbb S^1)$, the space of square-integrable functions on the circle, and that
the dynamical system $F$ in Eq.~\eqref{eq: Chafee-Infante} has a connected trapping region.
For $u\in X$, $F^t(u)\in C^\infty(\mathbb S^1)$ for $t>0$.
Hence, every function $u$ belonging to the global attractor is $C^\infty$.
Robinson also discusses more general cases where $f$ satisfies some general conditions, including $x f(x)<0$ for all sufficiently large $|x|$.

\medskip
We make the restriction  $0 <\lambda < 1$ because, in this range, the graph of Eq.~(\ref{eq: Chafee-Infante}) coincides with its ODE counterpart, Eq.~\eqref{eq: ODE}.
\black
It has exactly three nodes and each node consists of a constant function $u$.
The two nodes $u\equiv+1,-1$ are bottom nodes; 
the node $u\equiv0$ is a top node and the graph has two edges: one from 0 to 1 and one from 0 to -1.

\smallskip
\begin{theorem}[\bf Global attractor and graph]
    \label{thm:CI}
    Assume $0<\lambda<1$. 
    Then the global attractor of the Chafee-Infante PDE with $2\pi$-periodic boundary conditions consists of constant functions, where the constant is in $[-1,1]$. 
    In particular, the nodes consist of three constant functions $u(t,x)=1,-1,0$, where 0 is a top node and $1$ and $-1$ are bottom nodes (see page~\pageref{node} for definitions).
\end{theorem}
\smallskip

Each solution of the ODE~\eqref{eq: ODE} corresponds to a trajectory of the PDE for which $u_{xx}\equiv0.$
The global attractor of the ODE is the interval $[-1,1]$.
Hence, there is a one-to-one dynamics-preserving correspondence between the global attractor of the ODE and the global attractor of the PDE.

\medskip
We will use below the following Poincar\'e inequality, which is stated in more generality in the Appendix \blue\ref{sec:P}\black:
\beqn
\label{eq:Poincare}
\|u_x\|_{L^2(\mathbb S^1)}\leq\|u_{xx}\|_{L^2(\mathbb S^1)}\text{ for each }u\in H^2(\mathbb S^1).
\eeqn

\noindent
The proof involves the following three Lyapunov functions:
$$
V(u) \bydef \frac{1}{2}\|u_x\|^2_{L^2(\mathbb S^1)};
$$
$$
W_1(u) \bydef \max_x u(x);
$$
$$
W_2(u) \bydef \max_x -u(x).
$$
We also use the notation
$$
\frac{d^+}{dt}w(t)=\limsup_{\delta\to0^+}\frac{w(t+\delta)-w(t)}{\delta}.
$$
Recall that $f(u)=u(1-u^2)$.
\begin{lemma}
    \label{lemma: W}
    Under the hypotheses of Theorem~\ref{thm:CI},
    $
    \frac{d^+}{dt}W_j(u(t)) < 0
    $
    when $W_j(u(t))>1$, $j=1,2$.
\end{lemma}
We discuss the case $j=1$ where  $u(x)>1$, and the case for $j=2$ follows analogously.

\noindent
{\em Heuristic proof.}
For the moment, assume that $\frac{d}{dt} W_1(u(t))$ exists.
Whenever $u(t,x)=W_1(u)>1$, then $u_x(t,x)=0$ and $u_{xx}(t,x)\leq0$ and
$f(u(t,x))<0$. 
Hence 
$$\frac{d}{dt}W_1(u(t))\le f(u(t,x))<0.$$

\noindent
{\em Proof.}
Choose $t_n\bydef t+\delta_n\to t$ 
such that 

$$
\frac{W_1(u(t+\delta_n))-W_1(u(t))}{\delta_n}\to\frac{d^+}{dt}W_1(u(t)).
$$

Furthermore $x_n$ can be chosen to be a max of $u(t_n)$
so that 
$$W_1(u(t+\delta_n))= u(t+\delta_n, x_n).$$

We can assume that $x_n\to$ some $x_m$, where $x_m$ is a maximum of $u(t)$.

Hence 
$$
\frac{W_1(u(t+\delta_n))-W_1(u(t))}{\delta_n}
=
\underbrace{\frac{u(t+\delta_n,x_n)
-u(t,x_n)}{\delta_n}}_{\to u_t(t,x_m)\text{ as }n\to\infty}
+
\underbrace{
\frac{u(t,x_n)
-u(t,x_m)}{\delta_n}
}_{\leq0}
$$
so that 
$$
\frac{d^+}{dt}W_1(u(t))\leq u_t(t,x_m)
=\underbrace{u_{xx}(t,x_m)}_{\leq0}+\underbrace{f(u(t,x_m))}_{\leq f(W_1(u(t)))}<0.\qed
$$

\noindent{\em Proof of Theorem~\ref{thm:CI}.}

Now that we showed that each solution is uniformly bounded, we show that each solution converges to a constant function.
Indeed,
\smallskip
$$
\frac{d}{dt}V(u(t)) = 
$$
$$
=\langle u_x,\dot u_x\rangle_{L^2}
$$
(integrating by parts)
$$
=
-\langle u_{xx},\dot u\rangle_{L^2}
$$
\smallskip
$$
=
-\langle u_{xx},u_{xx}+\lambda f(u)\rangle_{L^2}
$$
\smallskip
$$
=
-\|u_{xx}\|^2_{L^2} - \lambda \langle u_{xx},f(u)\rangle_{L^2}
$$
(integrating by parts)
$$
=
-\|u_{xx}\|^2_{L^2} + \lambda \langle u_{x},\underbrace{f'(u)}_{1-3u^2}u_x\rangle_{L^2}
$$
\smallskip
$$
=
-\|u_{xx}\|^2_{L^2} - 3\lambda \langle u_x,u^2 u_x\rangle_{L^2} + \lambda\langle u_x,u_x\rangle_{L^2}
$$
\smallskip
$$
=
-\|u_{xx}\|^2_{L^2} - 3\lambda\|u u_x\|^2_{L^2}+ \lambda\|u_x\|^2_{L^2}.
$$
Therefore, by the Poincar\'e inequality (\ref{eq:Poincare}), 
\begin{align*}
    \frac{d}{dt}V(u(t)) \leq& -\|u_{xx}\|^2_{L^2} - 3\lambda \|uu_x\|^2_{L^2} + \lambda \|u_{xx}\|^2_{L^2} \\
    =& (\lambda-1)\|u_{xx}\|^2_{L^2}- 3\lambda \|uu_x\|^2_{L^2}\leq(\lambda-1)\|u_{xx}\|^2_{L^2},
\end{align*}
which is strictly negative for all $0<\lambda<1$ unless $u$ is a constant function of $x$ and then it is $0$.
\black
Hence, $\frac{d}{dt}V(u)\leq0$. 

For such $\lambda$, therefore, any solution $u(t)$ that is defined for all $t\geq0$ converges to a solution whose first space derivative is identically zero. Together with the boundary periodic conditions, 
this means that $u(t)$ must converge to a function constant in space. 
By Lemma~\ref{lemma: W}, this constant in absolute value cannot be larger than 1.

The global attractor includes only functions that are constant and, on the constant solutions, $u_{xx}=0$, and so the dynamics reduces to the ODE $\dot u=\lambda f(u)$. 
In this case, $[-1,1]$ is the global attractor of the ODE.
Hence, the global attractor of the PDE is the set of constant functions where the constant is in $[-1,1]$.
\qed

\medskip
More generally, consider the system of $n$ reaction-diffusion PDEs
\beqn
\label{eq:rd}
\partial_t u_k = D_k \Delta u_k + \lambda F_k(u_1,\dots,u_n),\,k=1,...,n,
\eeqn
where $D_k,\lambda$ are positive constants, $F=(F_1,\dots,F_n)$ is a smooth vector field  and the $u_k=u_k(t,x):\R\times\R^m\to\R$ are periodic in each space variable.
We conjecture that, for $\lambda$ small enough, if the ordinary differential equation $\dot x=F(x)$ has a compact connected global attractor, the graph of the ODE coincides with the graph of the PDE.
Furthermore, the global attractor of \eqref{eq:rd} consists solely of constant functions.

\section{Proofs}
\label{sec: proof}

In this section we will prove the main claims we made in Section~\ref{sec: chain graph}.
For the readers convenience, we repeat the statements of the propositions.

\subsection{\bf Global attractor.}
Recall that the global attractor $G$ of a dynamical system $F$, when it exists, is a ``maximal'' compact invariant set that attracts each compact set of the phase space $X$.
It is maximal in the sense that there is no strictly larger compact invariant set.

\medskip
\begin{proposition}[\bf Global attractor maximality]
\label{prop: glob}
The global attractor $G$ contains every compact invariant set of $F$ and is contained in every set that attracts all compact subsets of $X$.
\end{proposition}
{\bf Proof.}
Let $G'$ be a compact invariant set.
Then, for each $\eps>0$, there is $\tau>0$ such that $G'=F^t(G')\subset\Nbhd(G,\eps)$ for every $t\geq\tau$, since $G$ attracts every compact set.
Hence, $G'\subset G$.
Assume now that a set $A\subset X$ attracts each compact set. 
In particular, it attracts $G$ and so, by the same argument above,  $G\subset A$.\qed
\medskip

\subsection{\bf Global trapping regions.}
Recall that a global trapping region $Q$ for a dynamical system $F$ is a compact and forward invariant subset of $X$ that absorbs each compact subset of $X$.

The following result is our main motivation for introducing global trapping regions.

\medskip\noindent
{\bf Proposition 1 (Existence and connectedness of the global attractor).} {\em
Assume that $F$ has a connected global trapping region $Q$.
Then the global attractor exists and is connected -- and equals $\omega(Q)$.}

\medskip\noindent
{\bf Proof.}
Since $Q$ is closed and forward invariant, $\omega(Q)$ coincides with the nested intersection $\cap_{t\geq0}F^t(Q)$.
By assumption, the intersection is non-empty and compact, and it attracts $Q$.
Since $Q$ absorbs each compact subset of $X$, it follows that $\omega(Q)$ attracts each compact subset of $X$.
Since the intersection is nested, $\cap_{t\geq0}F^t(Q)=\cap_{t\geq s}F^t(Q)$ for any $s\geq0$, namely $\omega(Q)$ is invariant.
Hence, $\omega(Q)$ is indeed the global attractor of $F$.
The connectedness of $\omega(Q)$ follows from the following standard topological result: the nested intersection of connected compact sets is non-empty, compact and connected.
\qed

\medskip\noindent
{\bf Proposition 2.}
    {\em\ 
    Assume that a dynamical system $F$ has a global trapping region.
    Then the chain-recurrent set of $F$ is contained in the global attractor and the graph of $F$ is the same as the graph of the restriction of $F$ to its global attractor.
    In particular, all maximal chain-recurrent sets of $F$ are compact.
}

\begin{proof}
    We restrict $F$ to a global trapping region $Q$.
    Norton~\cite{Nor95} proved that, in the discrete-time case, the chain-recurrent set $R$ of a discrete dynamical system on a compact space is compact and invariant, and so are all of its nodes.
    By a result of Hurley~\cite{Hur95}, the chain-recurrent set of a continuous-time dynamical system is identical to the chain-recurrent set of its time-1 map, so Norton's result also holds for continuous-time systems.
    Hence, $R$ is contained in the global attractor of $F$.
    While Conley did not discuss graphs, he laid the framework by developing the theory of $\eps$-chains~\cite{Con72}.
    Moreover, Conley~\cite{Con78} proved that the chain-recurrent set of the restriction $F_R$ of $F$ to $R$ is $R$ itself, and each node of $F$ is a node of $F_R$ and vice versa.

    Now, suppose that there is an edge from $M$ to $N$ in $X$.
    We claim that $M$ is upstream from $N$ in the global attractor.
    To show this, we denote by $E$ the intersection for $\eps>0$ of the sets of all limit points of all $\eps$-chains from $M$ to $N$.
    
    We claim that the set $E$ is compact and invariant and, therefore, is contained in the global attractor. 
    Indeed, let $x\in E$.
    If $x\in M$ or $x\in N$, we know that there is a $y\in E$ with $F^1(y)=x$ because $M$ and $N$ are invariant. 
    If $x\in E\setminus(M\cup N)$, then it is the limit point of elements $x_n$, where $x_n$ is an element of a $1/n$-chain from $M$ to $N$ which is not the first of the chain.
    Since each $x_n$ is not the first element, then $x_n$ has a predecessor  $y_n$ in the $1/n$-chain, so that $d(F^1(y_n),x_n)<1/n$.
    By assumption, $x_n\to x$.
    Since $Q$ is compact, we can assume that $y_n\to y\in E$.
    Hence, by continuity, $d(F^1(y),x)\leq0$, namely $F^1(y)=x$.
    Moreover, $E$ is closed by construction and so is compact.

    Let now $\eps>0$, $x\in M$, $y\in N$.
    Let $E_\eps$ be the $\eps$-neighborhood of $E$.
    For $\eps>\eps^*>0$ sufficiently small, there is an $\eps^*$-chain $a_0,a_1,\dots,a_N$ from $x$ to $y$ lying totally in $E_\eps$.
    We can assume $\eps^*$ is sufficiently small that, if $p_1,p_2\in K$ and $d(p_1,p_2)\leq\eps^*$, then $d(F^1(p_1),F^1(p_2))<\eps$.
     We can construct a $3\eps$-chain from $x$ to $y$ in $E$ in the following way.
    Set $e_0=x$.
    Set $e_1$ to be a closest point of $E$ to $a_1$.
    Then
    \beq
    d(F^1(e_0),e_1)\leq d(F^1(e_0),a_1) + d(a_1,e_1)< 2\eps.
    \eeq
    For $k<N$, set $e_k$ to be a closest point of $E$ to $a_{k}$.
    Then 
    \beq
    d(F^1(e_{k-1}),e_k)\leq  d(F^1(e_{k-1}),F^1(a_{k-1})) + d(F^1(a_{k-1}),a_k)+d(a_k,e_k) < 3\eps.
    \eeq
    Finally, set $e_N=y$.
    Then
    \beq
     d(F^1(e_{N-1}),e_N)\leq  d(F^1(e_{N-1}),F^1(a_{N-1})) + d(F^1(a_{N-1}),a_N) < 2\eps.
    \eeq
    Hence, $N$ is downstream from $M$ in the global attractor.
\end{proof}
\blue
\medskip\noindent
{\bf Proposition 4.}
    {\em\ 
    Assume that a dynamical system $F$ has a global trapping region.
    Let $T>0$ and denote by $f$ its time-$T$ map.
    Then $f$ has the same chain-recurrent set, the same nodes, the same edges, and the same graph as $F$.
}

\medskip
To simplify notation, we prove the statement for $T=1$. 

\begin{proof}
    It is enough to prove that, if $x$ is upstream from node $N$ for $F$, then it is so also for $f$.
    Here we prove a simpler case, namely that if $x$ is chain-recurrent for $F$ it is so also for $f$.
    This case is a key lemma for the proposition.
    Of course all points in the node are chain-recurrent.
    We leave it to the reader  to complete the argument, as done in~\cite{DY25}.

    Let $\eps>0$. 
    By hypothesis, there is a (piecewise continuous) curve $z:[0,\tau]\to X$, $\tau\geq1$, such that $z(0)=z(\tau)=x$ and $z$ is $\eps$-close to the orbits of $F$.
    We want to prove that, for each $\eps$-close curve $z(t)$ there is an $\eps$-chain of $f$ between the same points.
    Since $\eps>0$ is arbitrary, that will prove the result.

    Denote by $\hat z$ the concatenation of $z$ with itself infinitely many times.
    So, $\hat z$
    is $\eps$-close to the orbits of $F$ and $\hat z(k\tau)=x$ for every natural integer $k$. 
    Moreover, for every integer $p>0$, the sequence $(x,f(x),\dots,f^p(x)) = (\hat z(0),\hat z(1),\dots,\hat z(p))$ is an $\eps$-chain for $f$ from $x$ to $\hat z(p)$.
    
    If $\tau$ is integer, then the claim follows at once.
    Assume that it is not.
    Hence, if $\tau$ is rational, there is an integer $k$ such that $k\tau$ is integer, so that $\hat z(k\tau)=x$, and so the claim follows.
    Finally, assume $\tau$ is irrational.
    There is a $\delta$-neighborhood of $\tau$ such that, if $|t-\tau|<\delta$, then $d(z(t),z(\tau))<\eps$.
    Furthermore, by definition of the concatenation, if $|t-\tau|<\delta$, then $d(\hat z(t),\hat z(\tau))<\eps$.
    Moreover, there is an integer $p$ within $\delta$ of some $k\tau$, where $k$ is an integer. 
    Hence, $d(\hat z(p),x)<\eps$, so that the sequence $(x,f(x),\dots,f^{p-1}(x),x)$ is an $\eps$-chain from $x$ to itself.
\end{proof}
\black

\subsection{\bf Connectedness of the chain graph.}
We prove here that the chain graph is connected under the assumption that there exists a connected global attractor $G$.
Recall that $G$ is invariant, namely, for every $t\geq0$, $F^t(G)=G$.
Therefore, 
for each $p\in G$ and every $t\geq0$, there is a point $q$ such that $F^t(q)=p$.
It follows that each $p\in G$ has a backward trajectory in $G$ defined for all $t$.
This trajectory is not necessarily unique, and different choices of backward trajectories may yield different backward limit sets. 

\medskip\noindent
{\bf Theorem 1 (Connectedness of (chain) graphs).} {\em
    Assume a dynamical system has a global attractor that is connected. 
    Then its (chain) graph is connected.}

\medskip\noindent
{\bf Proof.}
Denote by $G$ the global attractor of a dynamical system $F$ and by $\Gamma$ its graph.
Let $\cM$ and $\cN$ be disjoint collections of nodes of $\Gamma$ such that each node belongs to either $\cM$ or $\cN$ and denote by $M$ and $N$ the sets of the corresponding points in $X$.
We need to prove that if $M$ and $N$ are disjoint, then there is an edge from some node of $\cM$ to some node of $\cN$ or vice versa.
In order to prove this, let us set
$C=G\setminus(M\cup N)$. 
Each $p\in C$ has two limit sets $\alpha(p)$ and $\omega(p)$, where $\alpha(p)$ may not be uniquely defined.

There are four possible cases:

\smallskip\noindent
{\bf Case 0:} there is a point $p\in C$ whose limit sets $\alpha(p)$ and $\omega(p)$ are in different collections for some choice of $\alpha(p)$.
In this case there is a directed edge from a node in either $\cM$ or $\cN$ to a node in the other set, and the proof is done. 

\smallskip\noindent
{\BF Case $M$:} for each $p\in C$, the limit sets of $p$ are in $M$.
That is, the limit set of $p$ is in $M$ and every backward trajectory through $p$ has its limit set in $M$.
Since $G$ is connected, for each $\eps>0$ there is a point $p\in C$ for which dist$(p,N)<\eps$. 
Hence, there is a weak edge between some node of $\cM$ and some node of $\cN$.

\smallskip\noindent
{\BF Case $N$:} for each $p\in C$, the limit sets of $p$ are in $N$.
The proof is directly analogous to the proof of case $M$.

\smallskip\noindent
{\bf Remaining Case:} 
For each $p\in C$, $\omega(p)$ and every $\alpha(p)$ are in the same node collection, that is, some $p\in C$ has all its limit sets in $M$ and some $p'\in C$ has all its limit sets in $N$.
We say $p$ is type $M$ if all its limit sets are in $M$ and is type $N$ if all its limit sets are in $N$. 
Let $M^+$ be the set of all type $M$ points union $\overline{M}$ and define similarly $N^+$.
Both are closed sets, they are disjoint, and their union is $G$. 
That is a contradiction, because $G$ cannot be the union of two closed sets that are disjoint from each other.
\qed

\blue
\section{Conjectures}
\label{sec: conj}

\noindent
{\bf Two kinds of chain points.}
Let $y$ be downstream from $x$.
We say that $y$ is a {\BF limit chain point for $x$} if, in the definition of downstream, it is possible to choose $T_\eps$ so that $T_\eps\to\infty$ as $\eps\to0$.
As an example, every point in the limit set $\omega(x)$ is a limit chain point of $x$.

We say that $y$ is a {\BF trajectory chain point for $x$} if, in the definition of downstream, it is possible to choose $T_\eps$ so that $T_\eps$ is bounded as $\eps\to0$.

For a point $x$, a point $y$ can be both and an example is any time $x$ and $y$ are in the same periodic orbit.
So, $T_\eps$ can be chosen to be bounded and it can be chosen to go to infinity.

\medskip\noindent
\begin{conjecture}
\label{c1}
Assume $F$ is a continuous-time system with a compact global attractor.
Let $f$ be the  discrete time system $f(x)\bydef F^T(x)$ for some $T>0$.
Then the discrete and continuous-time systems have the same chain-limit points.
\end{conjecture}

\medskip

A non-trivial example of this equality comes from the map on the unit circle $F^t(x) = e^{i2\pi t}x$.
Let $f$ be the time-1 map of $F$.
This is the identity map.
Every point is a fixed point but, for each $x,y$ in the unit circle, $y$ is a limit chain point of $x$, because, for each $\eps>0$, there is an $\eps$-chain from $x$ to $y$.

\medskip\noindent
{\bf Can a graph have a weak edge?}
We have tried hard to find an example of a map with a weak edge, without success. 

\medskip\noindent
\begin{conjecture} 
\label{c2}
Assume $F$ is a dynamical system with a compact global attractor.
Then each edge of the graph of $F$ is strong.
\end{conjecture}

\black
\section{Appendix: Poincar\'e Inequality}
\label{sec:P}
\begin{lemma}[Poincar\'e Inequality]
    The two inequalities below hold:
     \par 
    \begin{enumerate}
    \item Let $\Omega\subset\bR^n$ be a connected bounded open set. Then there exists a constant $C>0$ such that 
    $$
    \|v\|_{L^2(\Omega)}\leq C\|\nabla v\|_{L^2(\Omega)}
    $$
    for every $v\in H^1_0(\Omega)$.
    \item Denote by $\bar v$ the average of a function $v\in L^2(\mathbb S^n)$, where $\mathbb S^n=[0,2\pi]^n$ with periodic conditions on the boundary.
    Then
    $$
    \|v-\bar v\|_{L^2(\mathbb S^n)}\leq \|\nabla v\|_{L^2(\mathbb S^n)}
    $$
    for every $v\in H^1(\mathbb S^n)$.
    \end{enumerate}
\end{lemma}
\begin{proof}
    1. Let $f$ be a function with compact support on $\Omega$.
    Then
    $$
    \int_\Omega f^2(x) dx = \int_\Omega \left[\frac{\partial}{\partial{xj}}x^j\right] f^2(x) dx = 
    $$
    $$
    = - 2\int_\Omega x^j f(x) f_{x^j}(x) dx \leq C\int_\Omega |f(x)|\cdot|\nabla f(x)|dx.
    $$
    Hence
    $$
    \|f\|^2_{L^2}\leq C\|f\|_{L^2}\|\nabla f\|_{L^2}.
    $$
    Since functions with compact support are dense in $H^1_0(\Omega)$, the claim follows.

    2. 
    Let $(x^1,\dots,x^n)\in\bR^n (\bmod 2\pi$) be angle coordinates on $\mathbb S^n$.     Set
    $$
    v(x) = \bar v + \sum_{k\neq0}a_k e^{ik\cdot x},\;k\in\mathbb Z^n.
    $$
    Then 
    $$
    \nabla v(x) = \sum_{k\neq0}k a_k e^{ik\cdot x}
    $$
    and, letting $\|k\|=\sqrt{k_1^2+\dots+k_n^2}$, 
    $$
    \|v-\bar v\|^2_{L^2(\mathbb S^n)}=\sum_{k\neq0}|a_k|^2
    \leq\sum_{k\neq0}\|k\|^2|a_k|^2 = \|\nabla v\|^2_{L^2(\mathbb S^n)}.
    $$
\end{proof}

\section*{Acknowledgments}

We are grateful to Joseph Auslander, Carlos Rocha and the referees for comments that helped us improve the manuscript.
This work was supported by the National Science Foundation, Grant Nos. DMS-1832126 and DMS-23082245.

\bibliography{refs}  

\end{document}